\documentclass[12pt]{amsart}
\usepackage[lite]{amsrefs} 
\usepackage{amssymb,amsmath,enumitem}
\DeclareMathOperator{\capa}{cap}
\newcommand{\defeq}{\mathrel{\mathop:}=}
\newcommand{\K}{{\mathbf K}}
\newcommand{\G}{{\mathbf G}}
\newcommand{\W}{{\mathbf W}}
\newcommand{\R}{{\mathbb R}}
\newcommand{\M}{\mathcal{M}^{+}}
\newcommand{\Om}{{\Omega}}

\numberwithin{equation}{section}
\newtheorem{theorem}{Theorem}[section]
\newtheorem{lemma}[theorem]{Lemma}
\newtheorem{remark}[theorem]{Remark}
\newtheorem{cor}[theorem]{Corollary}

\newtheorem{defn}[theorem]{Definition}

\allowdisplaybreaks
\begin{document}
\title[Nonlinear potential estimates] 
{Nonlinear potential estimates for sublinear problems with applications to 
elliptic semilinear and quasilinear equations}
\author{ Igor E. Verbitsky}
\address{Department of Mathematics, University of Missouri, Columbia, MO  65211, USA}
\email{verbitskyi@missouri.edu}
%%%%%%%%%%%%%%%%%%%%%%%%%
%%%%%%%%%%%%%%%%%%%%%%%%%
\begin{abstract}
	We give a survey of nonlinear potential estimates and their applications  
	obtained recently  for positive solutions to sublinear problems of the type 
		\[ 
		u = \mathbf{G}(\sigma u^q) + f  \quad \textrm{in} \,\, \Omega,
		\]
where $0 < q < 1$,  $\sigma\ge 0$ is a Radon measure in $\Omega$, $ f \ge 0$ is a measurable function, and 
$\mathbf{G}$ is a linear  integral operator with positive  
 kernel $G$ on $\Omega\times\Omega$. For  quasi-metric (or quasi-metrically modifiable) kernels $G$, these bilateral pointwise estimates  yield existence  criteria and uniqueness of solutions $u \in L^q_{{\rm loc}} (\Omega, \sigma)$. 
 
Applications  are considered to  semilinear elliptic equations involving the 
(fractional)  Laplacian,  
		\[ 
		(-\Delta)^{\frac{\alpha}{2}} u = \sigma  u^q + \mu \quad \textrm{in} \,\, \Omega, \qquad u=0 \, \, \textrm{in} \,\, \Omega^c. 
		\]
 Here $0<q<1$,  $\mu, \sigma \ge 0$ 
are   Radon measures, and $\Omega$ is a bounded uniform domain  in  $\R^n$, if $0 < \alpha \le 2$, 
 or  the entire space $\R^n$, a ball or half-space, if $0 < \alpha <n$. 
  
 Analogues of these results are presented for  elliptic equations involving the $p$-Laplace operator on the entire space $\R^n$,  
 \[ -\Delta_p u = \sigma  u^q + \mu \quad \textrm{in} \,\, \R^n, \qquad \liminf_{x\to \infty} u(x)=0, \]
 where $0<q<p-1$, and $\mu, \sigma \ge 0$ are   Radon measures.  More general 
 quasilinear equations with $\mathcal{A}$-Laplace operators ${\rm  div} \mathcal{A}(x, \nabla u)$ 
 in place of $\Delta_p$  are covered as well.
	\end{abstract}
\maketitle

\vfill 
\eject

%\tableofcontents

\section{Introduction}

\subsection{Sublinear problems} We give a survey of recent  developments 
in the theory of  sublinear problems based on nonlinear potential methods    (\cite{CV1}, \cite{PV2}, \cite{QV2}, \cite{Ver1}, \cite{Ver2}). In particular, we focus on  
	bilateral  pointwise estimates of solutions, which yield existence criteria and uniqueness 
	 results, for sublinear  integral equations  of the type 
	\begin{equation}\label{sublin-eq-f} 
			u =  \G(u^q d \sigma) + f, \quad  0 < u < \infty \quad d \sigma\textrm{-a.e.} \quad \textrm{in} \, \, \Omega,  
				\end{equation}
	where $0 < q < 1$, $f \ge 0$ is a Borel measurable function, and $\sigma\in \mathcal{M}^+(\Omega)$, the cone   
	of locally finite Radon measures in $\Omega$. (The class of finite Radon measures, with $\Vert \mu\Vert\defeq\mu (\Omega)<\infty$, is denoted by 
	$\mathcal{M}_{b}^+(\Omega)$.)

Here $\Omega$ is a locally compact Hausdorff space with countable base, and $\G$ is a linear integral operator with nonnegative kernel $G$ on $\Omega$,   
	$$\G^\sigma f (x)=\G (f \, d \sigma)(x) \defeq \int_{\Omega} G(x, y) \, f(y) d \sigma(y), \quad x \in \Omega,$$ 
	where $f \in L^1_{{\rm loc}} (\Omega, \sigma)$. If $f \equiv 1$,  we use the notation  
	\[
	\G \sigma(x) \defeq \int_{\Omega} G(x, y) \, d \sigma(y), \quad x \in \Omega, 
	\]
	for the linear $G$-potential of $\sigma\in \mathcal{M}^+(\Omega)$. 

	Throughout this paper, we use the following conventions 
	imposed on the \textit{kernels} $G$. 

	\begin{defn} A  kernel $G$ on $\Omega$ is understood to be  a  lower semicontinuous function  $G\colon \Omega\times \Omega\rightarrow [0, +\infty]$.  A  kernel $G$ is said to be  positive if  $G(x, y)>0$ for all $x, y \in \Omega$.
	\end{defn}

	We observe that, for equations \eqref{sublin-eq-f},  the sublinear case $0<q<1$ 
	is quite different  both from the linear case 
	$q=1$  and the superlinear case $q>1$, treated, for instance,  in \cite{FNV} and \cite{KV}, respectively. 	
	In particular, when $0<q<1$,  no ``smallness'' assumptions on $\sigma$ are needed 
	in order for a nontrivial solution $u$ to exist. Moreover, in contrast to the case $q\ge 1$, nontrivial 
	solutions to the homogeneous equations  ($f=0$) are treated similarly to non-homogeneous equations ($f\not=0$).

	For  \textit{quasi-metric} kernels, or more general quasi-metrically modifiable  kernels 
	discussed below  (see also \cite{FNV}, \cite{H}, \cite{KV}), we obtain matching bilateral pointwise estimates of solutions $u$ to \eqref{sublin-eq-f}.

\begin{defn}\label{q-m-def} A positive kernel $G$  on $\Omega$ 
	is said to be quasi-metric,  with quasi-metric constant $\kappa\ge \frac{1}{2}$, if $G$ is symmetric,  i.e., 
	$G(x,y)=G(y, x)$ for all  $x, y \in \Omega$, 
	and 
	 $d(x, y)\defeq \frac{1}{G(x, y)}$ satisfies  the quasi-triangle inequality   
		\begin{equation}\label{q-m-tr}
			d(x, y) \le \kappa [d(x, z) + d(z, y)], \qquad \forall \, x, y, z \in \Omega.
	\end{equation}
		\end{defn}
		
		Important examples of quasi-metric kernels include  Riesz kernels 
$I_\alpha (x, y)=c(\alpha, n) |x-y|^{\alpha -n}$  ($0<\alpha<n$), i.e., 
Green's kernels of the 
fractional Laplacian $(-\Delta)^{\frac{\alpha}{2}}$, on the 
entire Euclidean space $\mathbb{R}^n$ ($n \ge 1$), as well as Green's kernels  
of the Laplace--Beltrami operator on complete, non-compact Riemannian manifolds $M$ with nonnegative Ricci curvature (see, e.g.,  \cite{GSV}).

		The restriction that $G$ is symmetric in Definition \ref{q-m-def} 
		can be relaxed.

		\begin{defn}  A kernel $G$ on $\Omega$ is said to be quasi-symmetric {\rm (QS)}, with quasi-symmetry constant $\mathfrak{a}\ge 1$, 
	if 
		\begin{equation}\label{def-qs}
		\mathfrak{a}^{-1} G(y,x) \le G(x,y) \le \mathfrak{a} \, G(y,x), \quad \forall x, y \in \Omega. 
			\end{equation}
	\end{defn}
	
	The quasi-symmetry condition \eqref{def-qs} is often used below in combination with the \textit{weak maximum principle}.
	
	\begin{defn}	A kernel $G$ on $\Omega$ 
	satisfies the weak maximum principle {\rm (WMP)}, with constant $\mathfrak{b}\ge 1$, if 
			\begin{equation}\label{def-wmp}
			\G\mu (x) \le 1, \quad \forall x \in S_\mu \Longrightarrow \G\mu (x) \le \mathfrak{b},   \quad \forall x \in \Omega, 
				\end{equation}
		for any $\mu \in \mathcal{M}^+(\Omega)$,  where $S_\mu $ denotes the closed support of $\mu$. 
		\end{defn}

		When $\mathfrak{b}=1$, the kernel $G$ is said to satisfy the (Frostman) maximum principle.
			\smallskip

	\begin{remark}\label{rem-qm} Quasi-metric kernels are known to satisfy the {\rm (WMP)} with constant $\mathfrak{b} =2 \kappa$ 
(\cite{Ver2}*{Lemma 2.1}). 
\end{remark}

Some of the results for sublinear problems \eqref{sublin-eq-f}  hold for quasi-metric kernels. In particular, bilateral pointwise estimates of solutions to \eqref{sublin-eq-f}  are given 
	in terms of the linear potentials $\G \sigma$ and $\G (f \, d \sigma)$, as well as  certain ``intrinsic'' nonlinear potentials 
	$\K \sigma$ defined  below, which depend on $q\in (0, 1)$.  
	 
		Nonlinear potential  estimates lead to  the existence criteria for 
		\textit{all} solutions (possibly unbounded) to \eqref{sublin-eq-f}. They complement earlier results 
		 on the existence of positive solutions $u \in L^q(\Omega, \sigma)$ (globally)  in the homogeneous case $f=0$, which were 
		based on a sublinear version of Schur's lemma  for (QS)\&(WMP) kernels $G$ obtained in \cite{QV2}. 	
		
Bilateral pointwise estimates  
		also yield uniqueness of solutions 
		to \eqref{sublin-eq-f}  in the sublinear case ($0<q<1$). We observe 
		that  the uniqueness property may fail when $q\ge 1$.

More generally, we consider  nonlinear potential methods for 
 \textit{quasi-metrically modifiable} kernels $G$. 	 
\begin{defn} A positive kernel $G$  is said to be quasi-metrically modifiable,  
with modifier
		 $m\in C(\Omega)$, $m>0$, if the modified kernel 
		\begin{equation}
		\label{mod-ker}
		\widetilde{G}(x, y) \defeq \frac{G(x, y)}{m(x) \, m(y)}, \qquad x, y \in \Omega,
		\end{equation}
is quasi-metric, with quasi-metric constant $\tilde{\kappa}$. 
\end{defn}

A typical modifier for $G$   is given by 
	\begin{equation}
		\label{typ-mod}
g(x) = \min \{1, G(x, x_0)\}, \qquad x \in \Omega,
	\end{equation}
	where $x_0 \in \Omega$ is a fixed pole, provided $g\in C(\Omega)$.  

 Quasi-metrically modifiable kernels 	
 have numerous applications to semi-linear elliptic  PDE in domains 
 $\Omega$ with a positive Green's function $G$,  discussed in the next section 
 (see also  \cite{FNV}, \cite{FV1}, \cite{FV2},  \cite{GV}, \cite{H}, \cite{QV2}).

	We will discuss sharp lower estimates of (super) solutions, together with matching upper  estimates of (sub) solutions, 
	  to equation \eqref{sublin-eq-f},  
		 for quasi-metric, or quasi-metrically modifiable kernels, obtained in \cite{Ver2}*{Theorem 1.2}. As we will see below, lower estimates actually hold 
		 for  (QS)\&(WMP)   
	 kernels. 

\subsection{Semilinear elliptic equations}\label{subsec.1-2} Main applications of the results obtained for \eqref{sublin-eq-f} 
are concerned with sublinear elliptic equations of the type 
		\begin{equation}
			\label{frac_lap_eqn}
			\begin{cases}
				(-\Delta)^{\frac{\alpha}{2}} u =\sigma u^q + \mu, \,   & u>0 \, \text{ in } \Omega, \\
				u = 0 & \text{ in $\Omega^c$},
			\end{cases} 
		\end{equation}
		where $0<q<1$, $0<\alpha<n$, and $\mu, \sigma \in \mathcal{M}^+(\Omega)$,  in  a wide class of domains 
$\Omega$ in $\R^n$, or a Riemannian manifold, with Green's function $G$.

		If  $(-\Delta)^{\frac{\alpha}{2}}$  has a positive Green's function $G$ in  $\Omega$, then applying  Green's operator $\G$ to both sides, we obtain an equivalent problem where 
		solutions $u$ satisfy the integral equation 
		\eqref{sublin-eq-f} 
		with $f= \G \mu$. 
	If $\alpha=2$, 
	 such solutions $u$ to \eqref{frac_lap_eqn} 
	 in  bounded $C^2$-domains $\Omega$ 
	coincide with the so-called \textit{very weak} solutions  (see  \cite{MV}).

		Semilinear elliptic equations of this type have been  extensively studied, especially in the 
		classical case $\alpha=2$,  in bounded smooth domains $\Omega$ and on 
		$\R^n$,	 for bounded 
		solutions $u$, under substantial restrictions on the coefficients and data (see \cite{BoOr},   \cite{BO}, 
		\cite{Kr}*{Sec. 7.2.6}, 
		 and the literature cited there).  
		 
		 On the entire space $\R^n$, sharp existence and uniqueness 
		results were obtained by Brezis and Kamin \cite{BK} for \textit{bounded}  solutions $u>0$ to the equation $-\Delta u = \sigma u^q$. The proof of the uniqueness property given in \cite{BK} under the 
		assumption  $\displaystyle{\liminf_{x\to \infty} u(x)=0}$ is especially delicate. (Several simpler proofs are given in  \cite{BK} under the more restrictive condition 
		$\displaystyle{\lim_{x\to \infty} u(x)=0}$.) Pointwise estimates of bounded entire solutions given in \cite{BK} have a gap between the lower  and upper bounds.

		 Matching bilateral pointwise estimates were given recently in   \cite{Ver1}, \cite{Ver2} for all solutions $u \in L^q_{{\rm loc}} (\Omega, \sigma)$ to \eqref{frac_lap_eqn}, with arbitrary $\mu, \sigma \in \mathcal{M}^{+}(\Omega)$. As a consequence, the uniqueness problem was solved, and   sharp 
		existence criteria were given for such solutions, 
		in a certain class of bounded domains $\Omega\subset \mathbb{R}^n$ for $0<\alpha\le 2$, and on the entire space $\mathbb{R}^n$ 
		for $0<\alpha<n$, as well as on complete, non-compact Riemannian manifolds $M$ with nonnegative Ricci curvature (see \cite{GSV}).

	More precisely,  for $(-\Delta)^{\frac{\alpha}{2}}$, with $0<\alpha \le 2$ ($\alpha<n$)   
	in bounded uniform domains $\Omega\subset {\mathbb R}^n$, along with 
	 $\alpha=n=2$ in finitely connected domains $\Omega\subset {\mathbb R}^2$,  
	 Green's kernels 
	  are known to be quasi-metrically modifiable (see \cite{H}, and the literature cited there).  Hence, the general results  for \eqref{sublin-eq-f}
are applicable to all solutions  of \eqref{frac_lap_eqn} in these cases.

	When $0 < \alpha < n$, we can treat equations \eqref{frac_lap_eqn} for ``nice''  domains $\Omega$, such as the balls or half-spaces, where Green's kernel of $(-\Delta)^{\frac{\alpha}{2}}$ is known to be quasi-metrically modifiable  
	(see \cite{FNV}). 
	
	On the entire space $\Omega = {\mathbb R}^n$,  the Green kernel, i.e., the Newtonian kernel if $\alpha=2$, $n \ge 3$, and Riesz kernels of order $\alpha$ if $0<\alpha <n$, are quasi-metric. Equations  \eqref{frac_lap_eqn} in this case were 
	treated earlier in \cite{CV1} (existence  and bilateral pointwise estimates 
	for \textit{minimal} solutions).  
	 More complete results, including bilateral pointwise estimates 
	for \textit{all} solutions, and consequently 
	uniqueness of solutions,  
	were obtained subsequently in  \cite{Ver1} on ${\mathbb R}^n$, and in \cite{Ver2}  
	 for general quasi-metrically modifiable Green's kernels  
	 on $\Omega$.

\subsection{Quasilinear elliptic equations}\label{subsec.1-3} 	We will also present analogous nonlinear potential estimates, as  
well as existence and uniqueness theorems, for quasilinear 
elliptic equations involving  the $p$-Laplace operator, with lower order source terms, 
\begin{equation}\label{p_lap}
				-\Delta_p u =\sigma u^q + \mu, \quad   u>0 \, \text{ in } \, \, \R^n, \qquad \liminf_{x \to \infty} u = 0,
				\end{equation}
in the sub-natural growth case $0<q<p-1$. 

Bilateral estimates of all entire 
$p$-superharmonic solutions  (or, equivalently, local renormalized solutions)
were obtained  in \cite{Ver1}. They involve Havin--Maz'ya--Wolff potentials $\W\sigma$, $\W\mu$, and intrinsic 
nonlinear potentials $\K \sigma$, discussed in Sec. \ref{Sec-2}  and Sec. \ref{Sec-6} below. 

These pointwise estimates were used very recently to establish existence and uniqueness 
of the so-called \textit{reachable}  $p$-superharmonic solutions, in joint work with 
Nguyen Cong Phuc \cite{PV2}. 

More general  quasilinear  equations  of the type \eqref{p_lap},  
with $\mathcal{A}$-Laplace operators 
$\textrm{div} \mathcal{A}(x, \nabla u)$ in place of $\Delta_p $, under standard 
structural assumptions of order $p$ on  $\mathcal{A}(x, \xi)$ (see, e.g., \cite{HKM}),  and sub-natural growth terms, will be  discussed as well.

\subsection{A brief outline  of the paper}\label{subsec.1-4} 
Sec. \ref{Sec-2} contains some preliminary notions and basic concepts used throughout the paper. We first 
	consider certain weighted norm inequalities  of $(1, q)$-type for linear integral 
	operators $\G$ in $\Omega$. We then define intrinsic nonlinear 
	potentials $\K \sigma$,  
	using  localized versions of the $(1, q)$-type weighted norm  inequalities. We also discuss the precise definitions of sub- and 
super-solutions, and introduce the notion of the Wiener capacity 
for general kernels $G$.

In Sec. \ref{Sec-3}, we state the main results on  bilateral pointwise 
estimates, along with the existence and uniqueness results, for solutions to 
sublinear equations \eqref{sublin-eq-f}. We consider integral operators $\G$ 
with 
quasi-metric and  quasi-metrically modifiable kernels $G$, first 
 with data $f =\G \mu$, and then with arbitrary data $f \ge 0$. 

 In Sec. \ref{Sec-4}, we consider applications 
to semilinear elliptic problems of type \eqref{frac_lap_eqn} in uniform domains for $0<\alpha \le 2$, as well as the entire space, balls or half-spaces for $0<\alpha<n$. We also treat similar problems  involving linear uniformly  elliptic operators 
   with bounded measurable coefficients in non-tangentially accessible (NTA) domains,  or more general uniform domains with Ahlfors regular boundary.

In Sec. \ref{Sec-5}, we focus on the major steps in the proofs  
of the main theorems stated in  Sec. \ref{Sec-3} and provide relevant comments. In particular, we discuss 
 the key lemmas employed  in the proofs of the lower estimates of super-solutions for 
 (QS)\&(WMP) kernels $G$, and upper estimates of sub-solutions for 
quasi-metric and  quasi-metrically modifiable   kernels.

Finally, in Sec. \ref{Sec-6}, we are concerned  with analogous  nonlinear potential estimates  
and their applications 
for quasilinear equations \eqref{p_lap} involving the $p$-Laplacian, as well as more general $\mathcal{A}$-Laplace operators, 
in the case $0<q<p-1$. We  first treat  bilateral pointwise 
estimates for all  
 solutions, and then discuss the notion of a reachable solution, and present 
the corresponding existence and uniqueness results.  
 
 	\section{Preliminaries}\label{Sec-2}
	
\subsection{Weighted norm inequalities of $(1,q)$-type}\label{inequalities 1-q}

In this subsection, we discuss certain weighted norm inequalities studied 
in \cite{QV2}, along with their localized versions used extensively  below. 

We recall that throughout this paper, we use the notation $\mathcal{M}^{+}(\Omega)$ for \textit{locally finite} Radon measures in $\Omega$, and $\mathcal{M}^{+}_b(\Omega)$
for \textit{finite} Radon measures,  with $\Vert \nu \Vert\defeq \nu (\Omega)<\infty$ 
if $\nu \in \mathcal{M}^{+}_b(\Omega)$. All the kernels $G$  are assumed 
to be nonnegative lower semicontinuous functions defined on $\Omega\times \Omega$.

For  $\sigma \in \mathcal{M}^+(\Omega)$, $0 < q < 1$, and a kernel $G$  on $\Omega$, we consider weighted norm inequalities of $(1,q)$-type, 
\begin{equation}\label{main_ineq}
		\Vert \G \nu \Vert_{L^q(\Omega, \sigma)} \le C \, \Vert \nu \Vert, \qquad \forall \nu \in \mathcal{M}_b^+(\Omega). 
	\end{equation}

	 We denote by  $\varkappa=\varkappa(G, q, \sigma)$  the \textit{least constant} $C$ in \eqref{main_ineq}.

	  Clearly, \eqref{main_ineq} yields its $L^1$-version 
	 with $d \nu =f \, d \sigma$,  $f \in L^1(\Omega, \sigma)$, that is, 
	\begin{align}
		\label{main_ineq-1}
		\Vert \G^\sigma f \Vert_{L^q(\Omega, \sigma)} \le C \, \Vert f \Vert_{L^1(\Omega, \sigma)}, \qquad \forall f  \in L^1(\Omega, \sigma). 
	\end{align} 
	
\begin{remark} 	
	It follows from \cite{G}*{Lemma 3.I}  and  \cite{QV2}*{Theorem 1.1} that, conversely,    \eqref{main_ineq-1}$\Longrightarrow$\eqref{main_ineq}, 
	for  {\rm (QS)\&(WMP)} kernels $G$,  with a different constant $C$. However, it is easy to see that,  without the {\rm (WMP)} restriction, this implication may fail even for symmetric kernels $G$. 
		\end{remark}
	 
	 It is worth  observing that inequality \eqref{main_ineq-1} is the end-point case $p=1$ of the 
	$(p, q)$-type weighted norm inequality  
		\begin{equation}\label{p-r}
		\Vert \G^\sigma f \Vert_{L^q(\Omega, \sigma)} \le C \, \Vert f \Vert_{L^p(\Omega, \sigma)}, 
		\qquad \forall f \in L^p(\Omega, \sigma), 
	\end{equation}
	where $p \ge 1$ and $0<q<p$.  
	
	For $p>1$, $0<q<p$, inequality \eqref{p-r}
		was characterized recently in \cite{V1},  where it was shown that \eqref{p-r} holds,  for     kernels $G\ge 0$ that satisfy {\rm (QS)\&(WMP)} conditions, if and only if 
		\begin{equation}\label{p-r-crit}
		\int_\Omega (\G \sigma)^{\frac{q}{p-q}} d \sigma< \infty.
		\end{equation}
	
	In the more complicated case $p=1$, condition \eqref{p-r-crit} 
	is only necessary, but not sufficient, for \eqref{p-r} to hold. In \cite{QV2}*{Theorem 1.1}, it was  proved that, for (QS)\&(WMP)  kernels $G$, inequality \eqref{main_ineq}, or equivalently \eqref{main_ineq-1}, holds if and only if there exists a nontrivial super-solution $u\in L^q(\Omega, \sigma)$ of the homogeneous equation,  
	\[
			 u \ge \G(u^q d \sigma) \,\, d \sigma\textrm{-a.e.}  \quad  \text{ in $\Omega$}.  
		\]
	Moreover,  the least 
	constant $\varkappa$ in \eqref{main_ineq} satisfies the estimates 
	\begin{equation}\label{kappa-est} 
	 \Vert  u \Vert^{1-q}_{L^q (\Omega, \sigma)} \le \varkappa 
		 \le   C \, \Vert  u \Vert^{1-q}_{L^q (\Omega, \sigma)},  
	\end{equation}
	where  $C=C(q, \mathfrak{a}, \mathfrak{b})$ is a positive constant, and 
	$\mathfrak{a}$,   $\mathfrak{b}$ are
	 the constants in the 
	 conditions (QS), (WMP),  respectively. 
		
	This can be viewed as a sublinear  version of Schur's lemma  (see \cite{G}). 
	The proof is based on the notion of the  equilibrium 
	measure associated with the Wiener capacity for kernels $G$ discussed below (see \cite{Brelot},  \cite{Fug}). 
	
	In addition to estimates \eqref{kappa-est},  we have  (\cite{QV2}*{Theorem 1.2}),   
	\begin{equation}
		\label{lorentz}
	 C_1 \,   \Vert \G \sigma \Vert_{L^{\frac{q}{1-q}}(\Omega, \sigma)}  \le  \varkappa 
		 \le   C_2 \, 	 \Vert \G \sigma \Vert_{L^{\frac{q}{1-q}, q}(\Omega, \, \sigma)}, 
\end{equation}
	 where  $C_1=C_1(q, \mathfrak{b})$ and  
	 $C_2=C_2(q, \mathfrak{a}, \mathfrak{b})$ are positive constants. Here    $L^{r, q}(\Omega, \, \sigma)$ 
	 $(0<r<\infty, 0<q<\infty)$ 
	 stands for the Lorentz space on $\Omega$ with respect  to the measure 
	 $\sigma \in \M(\Omega)$.

\subsection{Intrinsic nonlinear potentials}\label{nonlin-pots} 
	
	Let $0<q<1$ and $\sigma \in \mathcal{M}^+(\Omega)$. Suppose $G$ is a kernel on $\Omega$. 
	In this section, we recall the definition of the intrinsic nonlinear 
	potential $\K \sigma$ given in \cite{Ver2}. Together 
	with the linear potential $\G \sigma$, it controls pointwise behavior of nontrivial solutions $u$ 
	to the homogeneous sublinear integral equation 
	\begin{equation}
			\label{sub-hom} 
			u= \G (u^q d \sigma), 
			\quad 0<u<+\infty \quad d \sigma\textrm{-a.e.},  \, \,  {\rm in} 
			\, \, \Omega.  
				\end{equation}
		 	 
	We define a ``ball'' $B= B(x, r)$ associated with $G$ by 
	 \begin{equation}\label{qm_ball}
	 B(x, r)\defeq\{ y\in \Omega\colon \, G(x, y)>1/r \}, \qquad  x\in \Omega, \, \, r>0.
	 \end{equation}
	 
	Notice that if $G$ is a quasi-metric kernel, then $B(x, r)$ is a \textit{quasi-metric ball} with respect to the quasi-metric $d=1/G$.

	\begin{remark}  By Fubini's theorem, $ \G\sigma$ can be represented in the form 
	\[
	 \G\sigma (x) = \int_0^\infty \frac{\sigma(B(x, r))}{r^2} \, d r, \qquad x \in \Omega. 
	\]
	 \end{remark}

	 Let   $d \sigma_{B} = \chi_{B} \, d \sigma$ be the restriction  of $\sigma$ to 
	 a ball $B$.  
	 We will need a localized version of inequality 
	 \eqref{main_ineq}, namely, 
	 \begin{align}
		\label{main_ineq_loc}
		\Vert \G \nu \Vert_{L^q(\Omega, \, \sigma_{B})} \le  C \, \Vert \nu \Vert, \qquad \forall \nu \in \mathcal{M}_{b}^+(\Omega), 
	\end{align}
	
	By $\varkappa(B)$ we denote the least constant  $C$ in \eqref{main_ineq_loc}. 
	We remark that  by \eqref{lorentz} with 
	$\sigma_B$ in place of $\sigma$, we deduce the following 
	 estimates of $\varkappa (B)$,  
	 \[
	 C_1 \,   \Vert \G \sigma_{B} \Vert_{L^{\frac{q}{1-q}}(\Omega, \, \sigma_{B})}  \le  \varkappa(B)
		 \le   C_2 \, 	 \Vert \G \sigma_{B} \Vert_{L^{\frac{q}{1-q}, q}(\Omega, \, \sigma_{B})}, 
	\]
	 where  $C_1=C_1(q, \mathfrak{b})$ and  
	 $C_2=C_2(q, \mathfrak{a}, \mathfrak{b})$ are positive constants. 
	 
	 The constants $\varkappa (B)$ with $B=B(x, r)$ are used to 
	  construct the nonlinear potential $ \K \sigma$,  intrinsic to 
	 the sublinear problem \eqref{sub-hom},  
	 \begin{equation}
			\label{nonlin_pot}
	 \K\sigma (x) \defeq \int_0^\infty \frac{\left[{\varkappa(B(x, r))}\right]^{\frac{q}{1-q}}}  {r^2}\,  d r, \qquad x \in \Omega. 
	 \end{equation}
	 
	 We remark that nonlinear potentials of this type were introduced for the first time in \cite{CV1} for Riesz kernels on $\Omega=\R^n$. In that case, $B=B(x, r)$ 
	 is a Euclidean ball of radius $r^{\frac{1}{n-\alpha}}$ centered at $x\in \R^n$. 
	
	Intrinsic nonlinear potentials $ \K\sigma$ 
	 resemble nonlinear potentials introduced originally by 
	 Havin and Maz'ya 
	 in \cite{MH}, but with $\sigma(B(x, r))$ used 
	 in place of $\varkappa(B(x, r))$. 
	 
	 More precisely, the Havin--Maz'ya--Wolff potential $ \W_{\alpha, p}$ 
	 on $\R^n$ (often called Wolff potential)  
	  is defined, for $0< \alpha<\frac{n}{p}$, $1<p<\infty$, by 
		 \begin{equation}
	 \label{nonlin-MH}
	 \W_{\alpha, p} \sigma (x) \defeq \int_0^\infty \frac{\left[{\sigma(B(x, \rho))}\right]^{\frac{1}{p-1}}}  {\rho^{\frac{n-\alpha p}{p-1}+1}}\, d \rho, \qquad x \in \R^n, 
	 \end{equation}
	 where $\sigma \in \M(\R^n)$, and 
	 $B=B(x, \rho)$ 
	 is a Euclidean ball in $\R^n$ of radius $\rho$ centered at $x$. 
	 
	 Notice that in the special case $p=2$, the potential $ \W_{\alpha, 2}$ coincides, up to 
	 a constant multiple, with the linear Riesz potential $\mathbf{I}_{2 \alpha}$ 
	 with kernel $I_{2\alpha}$  on $\R^n$.
	 
	  Nonlinear potentials  $ \W_{\alpha, p}$ 
	    were used subsequently by Hedberg and Wolff 
	 \cite{HW} in relation to the spectral synthesis problem for Sobolev spaces 
	 (see  \cite{AH},  and the literature cited there). In the special case $\alpha=1$, they 
	 are  fundamental to the theory of quasilinear elliptic equations 
	 of $p$-Laplace type  \cite{KM} (see 
	 also \cite{HKM}, \cite{KuMi}, \cite{Maz}, 
	 and Sec. \ref{Sec-6} below). 
	 
	 \subsection{Sub- and super-solutions}\label{sub-sup}

Let $\mu, \sigma \in \mathcal{M}^{+}(\Omega)$ and $0<q<1$.  
A Borel measurable function $u\colon \Omega \rightarrow [0, +\infty]$ is called a nontrivial \textit{super-solution}  associated with the equation 
\begin{equation}
			\label{sublin-sigma-mu} 
			u =  \G(u^q d \sigma) + \G \mu \quad  d\sigma\text{-a.e.} \, \, \text{in}  \, \, \Omega,  
			\end{equation}
if $u>0$ $d\sigma$-a.e., and 
\begin{equation}
			\label{super-sol-def} 
\G(u^q d \sigma) + \G \mu \le u <+\infty \quad d\sigma\text{-a.e.} \, \, \text{in}  \, \, \Omega.   
\end{equation}		
	 
	 A \textit{sub-solution} is defined similarly as a Borel measurable function $u\colon \Omega \rightarrow [0, +\infty]$ such that 
	 	 \begin{equation}
			\label{sub-sol} 
u \le \G(u^q d \sigma) + \G \mu<+\infty \quad d\sigma\text{-a.e.} \, \, \text{in}  \, \, \Omega.   
\end{equation}		

A  nontrivial \textit{solution} to \eqref{sublin-sigma-mu} is both a sub-solution and a nontrivial super-solution. 

If $u$ is a (super) solution, it is easy to see that actually $u\in L^q_{{\rm loc}} (\Omega, \sigma)$ (\cite{QV2}*{Lemma 2.2}).

 	 \subsection{The Wiener capacity}
	 
	 Let $G$ be a kernel on $\Omega$. For $\mu \in \M(\Omega)$, we set 
	\[
	\G^*\mu(y) \defeq \int_\Omega G(x, y) \, d \mu (x), \qquad y \in \Omega.
	\]
	Notice that the operator $\G^*$ is a formal adjoint of $\G$.

	Given a kernel $G$ on $\Omega$, a \textit{symmetrized} kernel $G^s$ is defined by
		\[ G^s(x,y) \defeq G(x,y) + G(y,x), \qquad x, y\in \Omega. \]
		Clearly,  $G^s$ is symmetric.  If $G$ is a (QS) kernel, then  $G^s$ is comparable to $G$: 
 \[ \left(1 + \mathfrak{a}^{-1}\right) \, G(x,y) \le G^s(x,y) \le (1 + \mathfrak{a}) \, G(x,y) , \qquad \forall x, y \in \Omega. \]
		
		The  kernel $G^s$ corresponds to the integral operator $ \G^s\defeq \G+\G^*$. 
	For a  (QS) kernel  $G$, the least constants  in the inequalities
			\begin{align*}
				& \Vert \G \nu \Vert_{L^q(\Omega, \sigma)} \le \varkappa \, \Vert \nu \Vert, \qquad \forall \nu \in \mathcal{M}^+(\Omega), \\
				& \Vert \G^s \nu \Vert_{L^q(\Omega, \sigma)} \le \varkappa_s \, \Vert \nu \Vert, \qquad \forall \nu \in \mathcal{M}^+(\Omega), 
			\end{align*}	
		are obviously equivalent: 
	$ \left(1 + \mathfrak{a}^{-1}\right) \, \varkappa \le \varkappa_s \le (1 + \mathfrak{a}) \varkappa$.

If  $G$ is a (QS) kernel, then there is a nontrivial  super-solution $u$, i.e., 		
	$\G(u^q d \sigma) +\G \mu \le u<\infty$ $d  \sigma$ a.e. 
			if and only if there is a nontrivial super-solution $u_s$ to the symmetrized version,   
			$\G^s(u^q_s d \sigma)+\G \mu  \le u_s<\infty$ $d  \sigma$ a.e.  
			This is easy to see using a scaled version $u_s= c_s \, u$ with an appropriate 
			constant $c_s>0$. A similar conclusion is true for sub-solutions. 
			
			Moreover, if $u_1$ is a sub-solution, and $u_2$ is a nontrivial 
			super-solution such that $u_1\le u_2$ then there exists a solution 
			$u$ such that $u_1\le u \le u_2$. 
			This means that if $G$ is a (QS) kernel, then 
			by passing to $G^s$, without loss of generality 
	we may assume that $G$ is  symmetric (see \cite{Ver2}).

	\begin{defn} 	
	The Wiener capacity $\capa (K)$ of a compact set $K \subset \Omega$ 
	is defined by 
		\begin{align} \label{capa}
			\capa (K) &\defeq \sup \, \{ \mu(K)\colon \mu \in \M(K), \, \,  \G^*\mu(y) \le 1,  \, \, \forall \, y \in S_\mu \}. 	
		\end{align}
		\end{defn} 

\begin{remark}
For positive symmetric kernels $G$, on any compact set $K \subset \Omega$ there exists an extremal measure $\mu$  (called equilibrium measure) such that $\mu(K)=\capa (K)<\infty$ in 
\eqref{capa} (see \cite{Brelot}, \cite{Fug}). Equilibrium measures play an important role in the proofs of the upper estimates of sub-solutions discussed below.
\end{remark} 	

 \begin{defn}    	A measure $\sigma\in \M(\Omega)$ is said to be \textit{absolutely continuous}  
		with respect to the Wiener capacity if $\sigma(K)=0$ whenever 
		 $\capa (K)=0$, for any compact set $K \subset \Omega$. 
		 \end{defn} 
	
	\begin{remark}\label{soln_abs_cont}
		Let $0<q< 1$ and $\sigma \in \M(\Omega)$. Let $G$ be a kernel on $\Omega$. If $u$ is a notrivial super-solution for $\G^*$ in place of $\G$, i.e.,  $u>0$ $d \sigma$-a.e. and  
$\G^*(u^q d \sigma) \le u<\infty$ $d\sigma$-a.e., then $\sigma$ is absolutely continuous with respect to the Wiener capacity  (see \cite{QV2}*{Lemma 4.2}).
	\end{remark}
	
	For (QS) kernels, the preceding remark is clearly true for nontrivial  super-solutions 
	$\G(u^q d \sigma) \le u<\infty$ $d\sigma$-a.e.

	\section{Main results for sublinear integral equations}\label{Sec-3}

	\subsection{Quasi-metric kernels}\label{main-q-m} We state our main theorem for quasi-metric kernels $G$.	
		
	\begin{theorem}\label{strong-thm}
		Let $\mu, \sigma \in \mathcal{M}^{+}(\Omega)$ ($\sigma\not=0$) and $0<q<1$.  Suppose $G$ is a  quasi-metric kernel on $\Omega$ 
		with quasi-metric constant $\kappa$. Then the following statements hold. 
		
		{(i)} Any nontrivial solution $u$ to equation \eqref{sublin-sigma-mu} 
					satisfies the bilateral pointwise estimates  
				\begin{equation}
				\label{sublin-low} 
		 c \,   [ (\G \sigma(x))^{\frac{1}{1-q}} + 	 \K\sigma (x)] +  \G \mu(x) \le 	u (x), 
		 \end{equation}	
		 \begin{equation} \label{sublin-up} 
		u(x)  \le  \,  C \, [ (\G \sigma(x))^{\frac{1}{1-q}} + 	 \K\sigma (x) +  \G \mu(x)],   
				\end{equation}	
		$d \sigma$-a.e. in $\Omega$, 
		where $c=c(q, \kappa)$, $C=C(q, \kappa)$ are positive constants. Moreover, such a solution $u$ is unique.

		{(ii)}  Estimate \eqref{sublin-low}  holds for any nontrivial super-solution $u$      
		at all $x \in \Omega$ such that  
		\begin{equation}
			\label{sublin-low-x} 
 \G(u^q d \sigma)(x) + \G \mu(x) \le u(x). 
\end{equation}		
		Similarly, 
		 \eqref{sublin-up}  holds for any sub-solution $u$   at all  $x \in \Omega$ such that 
		\begin{equation}
			\label{sublin-up-x} 
u(x) \le \G(u^q d \sigma)(x) + \G \mu(x). 
\end{equation}

		{(iii)} A nontrivial (super) solution $u$ to \eqref{sublin-sigma-mu}  exists if and only if 
		the following three conditions hold:
		\begin{align}
	& \int_a^\infty \frac{\sigma(B(x_0, r))}{r^2} \, d r < \infty,  \label{cond-a} \\ 
	& \int_a^\infty \frac{\left[{\varkappa(B(x_0, r))}\right]^{\frac{q}{1-q}}}  {r^2}\, d r<\infty,
	\label{cond-b} \\ & \int_a^\infty \frac{\mu(B(x_0, r))}{r^2} \, d r < \infty,  
				 \label{cond-c}
				 \end{align}	 
	for some (or, equivalently, all)  $x_0\in \Omega$ and $a>0$.	  Any nontrivial solution $u$ satisfies  \eqref{sublin-low}, \eqref{sublin-up} 
		 at all $x\in \Omega$ such that 
		\begin{equation}
			\label{sublin-eq-x} 
u(x) = \G(u^q d \sigma)(x) + \G \mu(x). 
\end{equation}	 
In particular, 	 \eqref{sublin-eq-x} holds $d \sigma$-a.e.
	\end{theorem}
	
	We notice that, as in the linear case \cite{FV2}, given 
	 a solution $u$ to \eqref{sublin-eq-x} defined $d \sigma$-a.e.,  
	 we can set  
	\begin{equation}
			\label{def-eq-x} 
	\tilde u(x)  \defeq \G(u^q d \sigma)(x) + \G \mu(x), \qquad \forall x \in \Omega. 
	\end{equation}	 
	 Then clearly $\tilde u = u$ $d \sigma$-a.e., and   
	\[
	\tilde u(x)  = \G(\tilde u^q d \sigma)(x) + \G \mu(x), \qquad \forall x \in \Omega.  
	\]
Hence,  the representative $\tilde u $  is a solution to \eqref{sublin-eq-x}  
 \textit{everywhere} in $\Omega$.

	\begin{remark}
	1. Under the assumptions of Theorem \ref{strong-thm}, conditions \eqref{cond-a}--\eqref{cond-c}, are  equivalent to   $\G \sigma <+\infty$, $\K \sigma <+\infty$, and $\G \mu < +\infty$ 
		 $d\sigma$-a.e. Other  equivalent conditions 
		 are given in Lemma \ref{equiv-lemma} and Corollary \ref{equiv-cor} below.

 2. An analogue of Theorem \ref{strong-thm} holds for  equation \eqref{sublin-eq-f} 
	with arbitrary $f\ge 0$ (see Theorem \ref{main-thm-f} below). One only needs to replace $\G \mu$ with $\G(f^q d \sigma)+ f$ in \eqref{sublin-low}, 
	\eqref{sublin-up}, and the corresponding  estimates for  sub- and 
	super-solutions. 
	The  term $\G(f^q d\sigma)$ is no longer needed in the special case $f=\G \mu$. 
\end{remark}

	\subsection{Quasi-metrically modifiable kernels}\label{main-q-q-m} 
	
	We next state our main theorem for quasi-metrically modifiable  kernels $G$.	In this case,  Theorem~\ref{strong-thm}  holds  
		with $\widetilde{G}$ in place of $G$, which leads to 
		matching  lower and upper global  estimates of solutions up to the boundary of $\Omega$.

	The modification procedure 
	is applicable to  Green's kernels $G$  for $(-\Delta)^{\frac{\alpha}{2}}$ in some   domains $\Omega\subset \R^n$,  in particular, balls or half-spaces, if  $0<\alpha <n$, or uniform domains discussed below 
	if $0<\alpha\le 2$.  In these cases, we use modifiers $g$ defined by  \eqref{typ-mod}.
	 For bounded $C^{1,1}$-domains $\Omega$ if $0<\alpha\le 2$, 
	as well as  balls or half-spaces if  $0<\alpha <n$, it is known that $g(x) \approx [\textrm{dist}(x, \Omega^c)]^{\frac{\alpha}{2}}$.

	 Suppose $G$  is a quasi-metrically modifiable kernel, with modifier $m$, associated with the quasi-metric 
	$\widetilde{d}=1/\widetilde{G}$. 		
	We denote by $\widetilde{B}(x, r)$  a quasi-metric ball 	
	\begin{equation}\label{def-tilde-mod} 
	\widetilde{B}(x, r) \defeq \left\{ y \in \Omega\colon \, \, \widetilde{G}(x,y)>1/r \right\}, \quad x \in \Omega, \, \, r>0. 
		\end{equation}

	Let $d \tilde \sigma=m^{1+q} d \sigma$.  For a Borel set $E \subseteq \Omega$, 
	 by $\tilde{\varkappa}(E)=\tilde{\varkappa}(E, \tilde \sigma)$ we denote the least constant in the inequality 
	\begin{equation}\label{def-kap-K} 
	 		\Vert \widetilde{\G} \nu \Vert_{L^q(\Omega, \widetilde{\sigma}_E)} \le \tilde{\varkappa}(E) \, \Vert \nu\Vert,
			\qquad \forall  
			\, \,  \nu \in \M_b(\Omega).
				 	\end{equation}

	Using the constants $\tilde{\varkappa}(\widetilde{B}(x, r))$, we construct the modified intrinsic potential $\widetilde{\K} \sigma$  defined by 
		\begin{equation}\label{def-K-mod} 
	\widetilde{\K} \sigma  (x) \defeq \int_0^\infty \frac{ [\tilde{\varkappa}(\widetilde{B}(x, r))]^{\frac{q}{1-q}} }{r^2} \, dr, \qquad x\in \Omega.
		\end{equation}

	\begin{theorem}\label{strong-thm-mod}
		Let $\mu, \sigma \in \mathcal{M}^{+}(\Omega)$ ($\sigma\not=0$) and $0<q<1$.  Suppose 
		$G$ is a  quasi-metrically modifiable kernel with modifier  $m$.		Then any nontrivial solution $u$ to equation \eqref{sublin-sigma-mu} is unique and 
					satisfies the bilateral pointwise estimates  
				\begin{equation}
				\label{sublin-low-mod} 
		  c \,   m \,  \left( \left[m^{-1}\G (m^q d\sigma)\right]^{\frac{1}{1-q}} + 
		   \widetilde{\K} \sigma\right) +   \, \G \mu\le 	u,  
		   \end{equation} 
		   \begin{equation}
		\label{sublin-up-mod} 
		u \le  \,  C \, m \,  \left( \left[ m^{-1}\G (m^q d\sigma)\right]^{\frac{1}{1-q}}  + 	  \widetilde{\K} \sigma\right) +   C \,  \G \mu,   
				     \end{equation} 	
		$d \sigma$-a.e. in $\Omega$, 
		where $c, C$ are positive constants which depend only on $q$ and the 
		quasi-metric constant $\tilde{\kappa}$ of the modified kernel $\widetilde{G}$. 	
		
		The lower bound 	\eqref{sublin-low-mod} holds for any nontrivial super-solution $u$, whereas the upper bound \eqref{sublin-up-mod}  holds for any sub-solution $u$. 
		\end{theorem}

	\begin{remark}\label{ref-a-mod}   1. Under the assumptions of Theorem \ref{strong-thm-mod}, a 
	 nontrivial (super) solution to \eqref{sublin-sigma-mu} exists  if and only if 
	 $\G (m^qd\sigma) <+\infty$, $\widetilde{\K} \sigma <+\infty$, and $\G \mu < +\infty$ 
		 $d\sigma$-a.e.

2. If $G$ is quasi-metrically modifiable with  modifier  $m=g$ given by \eqref{typ-mod}, 
	a nontrivial  (super) solution to \eqref{sublin-sigma-mu} exists if and only if  
	(Lemma \ref{equiv-lemma-mod} below),  
	\begin{equation} \label{tilde-kappa}
		\tilde{\varkappa} (\Omega)<\infty \quad \textrm{and}    \quad 	\int_\Omega g  \, d \mu<\infty.  \end{equation}
	\end{remark}

	\subsection{Equations with arbitrary data}\label{Sec-data}

We now state the main theorem 
for equation \eqref{sublin-eq-f} 
		with arbitrary  data $f \ge 0$ in place of $\G\mu$ (see \cite{Ver2}*{Theorem 6.1}).

%%%%%%%%%%%%%%%%%

\begin{theorem}\label{main-thm-f}
		  Suppose $0<q<1$, $G$ is a  quasi-metric kernel on $\Omega$,  
		   $\sigma \in \mathcal{M}^{+}(\Omega)$ ($\sigma\not=0$) and 
		  $f\ge 0$ is a Borel measurable function in $\Omega$. Then the following statements hold. 
		
		{(i)} Any nontrivial solution $u$ to equation \eqref{sublin-eq-f}  is unique and 
					satisfies the bilateral pointwise estimates  
				\begin{align}
				\label{lower-f} 
		 &c \,   \left[ (\G \sigma)^{\frac{1}{1-q}} + 	 \K\sigma +  \G (f^q d \sigma)\right] +f   \le u,\\ 
		& u \le  \,  C \, \left[ (\G \sigma)^{\frac{1}{1-q}} + 	 \K\sigma + \G (f^q d \sigma)\right]+f,   	
		\label{upper-f} 
\end{align}	
		$d \sigma$-a.e. in $\Omega$, 
		where $c=c(q, \kappa)$, $C=C(q, \kappa)$ are positive constants.

		{(ii)}  The lower estimate \eqref{lower-f}  holds for any nontrivial super-solution,      
		whereas  the upper estimate in \eqref{upper-f}   holds for any 
		sub-solution.

		{(iii)} A nontrivial solution $u$ to \eqref{sublin-eq-f}   exists if and only if 
	\begin{equation}	 \label{exist-f}
	\G \sigma <\infty, \, \, \K \sigma <\infty, \, \,  f <\infty, \, \,  \G(f^q d \sigma) <\infty, \quad d \sigma\textrm{-a.e.} \, \, {\rm in} \, \, \Omega.
\end{equation}	 
	\end{theorem}
	
\begin{remark}  1. In the special case $f=\G \mu$  
($\mu \in \M(\Omega)$) of Theorem \ref{main-thm-f}, 
the term $\G(f^q d \sigma)$ in  \eqref{upper-f}  may be dropped if, at the same time, $f$  is replaced with $C  f$. Moreover, the condition $\G(f^q d \sigma)<\infty$ 
$d \sigma$-a.e. in \eqref{exist-f} is redundant (see Theorem \ref{strong-thm} above). 
 
 2.  An analogue of  Theorem \ref{main-thm-f}  holds  for quasi-metrically modifiable kernels $G$. This gives an extension of Theorem \ref{strong-thm-mod} to solutions of equation  \eqref{sublin-eq-f}  with arbitrary data $f\ge 0$. The corresponding estimates of solutions   remain valid once we replace $\G\mu$ with $c \, \G(f^q d \sigma)+f$ in \eqref{sublin-low-mod}, and $C \, \G(f^q d \sigma)+f$ in \eqref{sublin-up-mod}, respectively. In the existence criteria  
discussed in Remark \ref{ref-a-mod}, it suffices  to replace $\G\mu$ with  $\G(f^q d \sigma)$. 
 \end{remark}

 \section{Semilinear elliptic problems}\label{Sec-4} 

 \subsection{Equations with the fractional Laplace operator}	 In the following definition of a \textit{uniform domain} (or, equivalently, an interior NTA domain), we rely on the notions of the interior corkscrew condition and the Harnack chain condition. We refer to 
   \cite{H} for related definitions, in metric spaces, along with 
	a discussion of quasi-metric properties, 3-G inequalities,   and the uniform boundary Harnack principle (see also \cite{Ai}).

 \begin{defn}
 A uniform domain $\Omega\subset \R^n$, $n \ge 2$,  
is a bounded domain which 
satisfies  the interior corkscrew condition and the Harnack chain condition. 
\end{defn}

	Notice that uniform domains are not necessarily regular in the sense of Wiener. Bounded Lipschitz and non-tangentially accessible (NTA) domains are examples of regular uniform domains.
	
	The next theorem (\cite{Ver2}*{Theorem 1.2}) is a direct consequence of Theorem \ref{strong-thm-mod} and the fact that Green's 
 function $G$ of  $(-\Delta)^{\frac{\alpha}{2}}$  in  a uniform domain   
 $\Omega$ for  $0<\alpha \le 2$  is quasi-metrically modifiable, with $m=g$ and quasi-metric constant $\tilde{\kappa}$ 
 which does not depend on the choice of $x_0\in \Omega$  
  (see \cite{An}, \cite{H}).

	\begin{theorem}\label{appl-unif} Suppose $\Omega\subset \R^n$, $n \ge 2$,  is a uniform domain. Suppose   $G$ is   Green's kernel of $(-\Delta)^{\frac{\alpha}{2}}$  in $\Omega$, where $0<\alpha \le 2$,   
	 $\alpha <n$. Define the modifier $m=g$  by 
 \eqref{typ-mod} with pole $x_0\in \Omega$. 
	
	  Let $0<q<1$, and let $\mu, \sigma \in \M(\Omega)$, and $d \tilde \sigma = 
	g^{1+q} d \sigma$.   
Then the following statements hold. 
 
  (i) Any nontrivial solution $u$ to equation \eqref{frac_lap_eqn} is unique and 
 satisfies estimates \eqref{sublin-low-mod}, \eqref{sublin-up-mod}  $d \sigma$-a.e., 
 and at all $x \in \Omega$ where \eqref{sublin-eq-x} holds. 
 
 (ii) Any nontrivial super-solution $u$ satisfies the lower bound \eqref{sublin-low-mod},  and any  sub-solution $u$ satisfies the upper bound \eqref{sublin-up-mod}.

(iii) A nontrivial (super) solution to  \eqref{frac_lap_eqn}  exists if and only if \eqref{tilde-kappa}
holds,  for some (or, equivalently, all) $x_0\in \Omega$.	 
\end{theorem}

	\begin{remark} Uniqueness of solutions to  sublinear problems of the type  \eqref{frac_lap_eqn} was 
	previously known only under heavy restrictions on solutions, coefficients, and data, for instance, for bounded solutions   \cite{BK}, \cite{BO}, or finite energy solutions \cite{SV}.	
\end{remark} 

For specific domains $\Omega\subseteq \R^n$, Theorem \ref{appl-unif} can be extended 
to the full range $0<\alpha<n$. 

\begin{theorem}\label{frac-special} Let $0<\alpha <n$. 
Suppose $\Omega$ is the entire space $\R^n$, or a ball, or half-space in $\R^n$. 
Then an analogue of Theorem \ref{appl-unif} holds. 
\end{theorem}

\begin{remark}
  In all the cases listed in Theorem \ref{frac-special}, Green's kernel $G$ of 
$(-\Delta)^{\frac{\alpha}{2}}$, for  $0<\alpha <n$,   is known to be either quasi-metric (Riesz kernel $I_\alpha$ of order $\alpha$ if $\Omega=\R^n$), 
or quasi-metrically modifiable ($\Omega$ is a ball, or half-space), with 
modifier $m(x)=g(x)\approx [{\rm dist} (x, \Omega^c)]^{\frac{\alpha}{2}}$    (see \cite{FNV} and the literature 
 cited there). 
  \end{remark}

 \subsection{Equations with uniformly elliptic operators}	A similar approach works also for Green's kernels $G$ of 
	uniformly elliptic, symmetric operators $\mathcal{L}$ in divergence form, 
	\begin{equation}
		\label{typ-L}
		 \mathcal{L} u=-\textrm{div} (A \nabla u), \quad A=(a_{ij}(x))_{i, j=1}^n, \quad a_{ij}=a_{ji} ,
			\end{equation}
with real-valued coefficients $a_{ij}\in L^\infty(\Omega)$, in place of 
$(-\Delta)^{\frac{\alpha}{2}}$. 
	 
 Then, as above,  Green's  kernel $G$ of $\mathcal{L}$ is known to be quasi-metrically modifiable, with modifier $m=g$  defined by 
 \eqref{typ-mod} with pole $x_0\in \Omega$, under 
 certain restrictions on $\Omega$  specified in the following theorem 
  (see \cite{FV2}, \cite{H}, \cite{Ver2}, and the literature cited there). 
  	
	\begin{theorem} An analogue of Theorem \ref{appl-unif},  for 
	operators $\mathcal{L}$ in place of
	$(-\Delta)^{\frac{\alpha}{2}}$,  holds  for linear uniformly  elliptic operators 
	$\mathcal{L}$ with bounded measurable coefficients 
 given  by \eqref{typ-L}, and modifiers $m=g$, 
 in NTA domains, as well as uniform domains with Ahlfors regular boundary. 
 \end{theorem}

	\section{Outlines of the proofs of the main theorems}\label{Sec-5}

%%%%%%%%

\subsection{Lower bounds for super-solutions: (QS)\&(WMP) kernels}\label{lower-qm} Let $G$ be a kernel on $\Omega$. 
We first consider nontrivial super-solutions to the 
homogeneous equation \eqref{sub-hom}, i.e., functions $u>0$ $d\sigma$-a.e. such that (see Sec. \ref{sub-sup})
\begin{equation}
			\label{super-sol-hom} 
\G(u^q d \sigma)  \le u <+\infty \quad d\sigma\text{-a.e.} \, \, \text{in}  \, \, \Omega.   
\end{equation}		
	
	We start with the following lower bound for super-solutions obtained in 
	\cite{GV}*{Theorem 1.3} for {\rm (WMP)} kernels. 
	\begin{lemma}
		\label{G-lemma-lower} Let $\sigma \in \M(\Omega)$ and $0<q<1$.  
		Suppose $G$ is a kernel on $\Omega$ which satisfies the {\rm (WMP)}  with constant 
		$\mathfrak{b}$ in \eqref{def-wmp}. Then any nontrivial super-solution $u$ 
		 satisfies the estimate
		\begin{equation}\label{G-lower-est} 
	 		u(x)\ge c \,  \left[\G\sigma(x)\right]^{\frac{1}{1-q}}, 
	 	\end{equation}
	where $c=(1-q)^{\frac{1}{1-q}} \mathfrak{b}^{-\frac{q}{1-q}}$, for all $x\in \Omega$ such that  $\G(u^q d \sigma)(x)\le u(x)$. In particular, \eqref{G-lower-est}  
	holds $d \sigma$-a.e.
	\end{lemma}
	
There is another lower estimate for super-solutions (\cite{Ver2}*{Lemma 3.2}),  which complements \eqref{G-lower-est} in a crucial way.   It holds 
for kernels $G$ which satisfy both the (WMP) and (QS) conditions. Using a symmetrized kernel, 
we may assume without loss of generality that $G$ is symmetric.
	\begin{lemma}
		\label{K-lemma-lower} Let $\sigma \in \M(\Omega)$ and $0<q<1$.  
		Suppose $G$ is a symmetric kernel on $\Omega$ which satisfies the {\rm (WMP)} with constant 
		$\mathfrak{b}$. Then any nontrivial super-solution $u$ satisfies the estimate
		\begin{equation}\label{K-lower-est} 
	 		u(x)\ge c \,  \K \sigma(x), 
	 	\end{equation}
	where $c=(1-q)^{\frac{1}{1-q}} \mathfrak{b}^{-\frac{q}{1-q}}$, for all $x\in \Omega$ such that $u(x)\ge \G(u^q d \sigma)(x)$. In particular, \eqref{K-lower-est} holds $d \sigma$-a.e.
	\end{lemma}

	For (QS)\&(WMP) kernels, the constant $c$ in  \eqref{K-lower-est}  will depend on $q$, $\mathfrak{b}$, and additionally the constant $\mathfrak{a}$ in 
	\eqref{def-qs}. Combining  \eqref{G-lower-est} with 
	 \eqref{K-lower-est}  and the trivial estimate $u \ge \G \mu$, 
we obtain the main lower estimate for 
any nontrivial super-solution $u$ to \eqref{sublin-sigma-mu}.
	
	\begin{cor} 
	\label{cor-lower-est} Let $\mu, \sigma \in \M(\Omega)$ and $0<q<1$.  
		Suppose $G$ is a ${\rm (QS)}\&{\rm (WMP)}$ kernel on $\Omega$ 
		with constants $\mathfrak{a}$, $\mathfrak{b}$  in \eqref{def-qs},  \eqref{def-wmp}, respectively. Then any nontrivial super-solution $u$ to \eqref{sublin-sigma-mu} satisfies the estimate
		\begin{equation}\label{K-lower-est-mu} 
	 		u(x)\ge c \,  [(\G \sigma(x))^{\frac{1}{1-q}}+\K \sigma(x)] + \G \mu(x), 
	 	\end{equation}
	where $c=c(q, \mathfrak{a}, \mathfrak{b})$, for all $x\in \Omega$ such that 
	\begin{equation}\label{x-lower-est-mu} 
	 	 u(x) \ge \G(u^q d \sigma)(x)+ \G \mu(x). 	
	 	\end{equation}
	In particular, \eqref{K-lower-est-mu} holds $d \sigma$-a.e.
	\end{cor}

	\subsection{Upper bounds for sub-solutions: quasi-metric kernels}\label{q-m kernels}
	
	In this section, we discuss the main steps involved in the proof of the upper estimates \eqref{sublin-up} 
	of sub-solutions associated with  equation  
	 \eqref{sublin-sigma-mu}, for quasi-metric kernels $G$. They match the lower estimates of super-solutions    obtained in Corollary \ref{cor-lower-est}.

	We start with the following key estimate (\cite{Ver2}*{Lemma 5.4}). 			
	\begin{lemma}
		\label{qm_lemma-upper}
		Let $G$ be a quasi-metric  kernel on $\Omega$ with 
		quasi-metric constant $\kappa$. Let $0<q<1$ and $\nu, \sigma \in  \M(\Omega)$. 
			Then, for all $x \in \Omega$,  
					\begin{equation}\label{upper-est-qm} 
	 	\G [(\G \nu)^q d \sigma] (x)\le C   \left (\G \nu(x) \right)^q \, 
		\left[ \G\sigma(x) 
			+ \left(\K \sigma (x)\right)^{1-q} \right], 
	 	\end{equation}
		where $C=(2 \kappa)^{q}$.
	\end{lemma}
	
The following lemma, which is  deduced from Lemma \ref{qm_lemma-upper}, 
yields the desired upper estimate for sub-solutions, but only  $d \sigma$-a.e. 
The remaining difficulty is to prove the upper estimate for points $x$ where 
possibly $u(x)=+\infty$; it is handled in 
Lemma \ref{qm-upper-K} below. 
	
	\begin{lemma}
		\label{qm_cor-upper}
		Let $G$ be a quasi-metric  kernel on $\Omega$ with 
		quasi-metric constant $\kappa$. 	Let $0<q<1$ and $ \mu, \sigma \in  \M(\Omega)$. 		Then any sub-solution $u \ge 0$ such that  $u\le \G(u^q d \sigma)+\G \mu<+\infty$ $d \sigma$-a.e.,  satisfies the estimate
		\begin{equation}\label{upper-est-qm-cor} 
	 		u(x)\le C \,   \left [ \left(\G\sigma(x)\right)^{\frac{1}{1-q} } 
			+ \K \sigma (x) + \G\mu(x)\right], 
	 	\end{equation}
		for all $x\in \Omega$ such that $u(x)\le \G(u^q d \sigma)(x) +\G\mu(x)<+\infty$, where 
		$C= (8 \kappa)^{\frac{q}{1-q}}$. In particular, \eqref{upper-est-qm-cor} holds $d \sigma$-a.e. 
	\end{lemma}

	In what follows,   we will repeatedly use the fact, mentioned in  Remark \ref{rem-qm} above,   
	   that a quasi-metric kernel with quasi-metric constant $\kappa$ obeys the (WMP) with constant 
	 $\mathfrak{b}=2\kappa$.

	The next lemma is used to deduce 	estimate	\eqref{upper-est-qm-cor}  for \textit{all} $x \in \Omega$ such that $u(x) \le \G(u^q d \sigma)(x)+\G \mu(x)$, including the case $u(x)=+\infty$.

	\begin{lemma}
		\label{qm-upper-K}
		Let $G$ be a quasi-metric  kernel on $\Omega$ with 
		quasi-metric constant $\kappa$. Let $0<q<1$ and $ \mu, \sigma \in  \M(\Omega)$. 
		 Then the function 
		 \begin{equation}\label{def-h} 
		  h(x)\defeq (\G\sigma(x))^{\frac{1}{1-q}} 
			+ \K \sigma (x) + \G \mu(x), \quad x \in \Omega,
		\end{equation} 
		 satisfies the estimate  
		 	\begin{equation}\label{upper-est-qm-K} 
	 	\G (h^q d \sigma) (x)\le C  \, h(x), \qquad \forall x \in \Omega, 
	 	\end{equation}
		 where $C$ is a constant which depends only on $q$ and $\kappa$.
	\end{lemma}
	
	Combining Lemma \ref{qm_cor-upper} and Lemma \ref{qm-upper-K}  
	yields the following corollary.
		
		\begin{cor}
		\label{qm_cor-upper-K}
		Let $G$ be a quasi-metric  kernel on $\Omega$ with 
		quasi-metric constant $\kappa$. 		Let $0<q<1$ and $ \mu, \sigma \in  \M(\Omega)$. 	Then every sub-solution $u$ for which $u\le \G(u^q d \sigma) + \G \mu<+\infty$ $d \sigma$-a.e.  satisfies the estimate
		\begin{equation}\label{upper-est-qm-cor-K} 
	 		u(x)\le (8 \kappa)^{\frac{q}{1-q}} \,   \left [ \left(\G\sigma(x)\right)^{\frac{1}{1-q} } 
			+ \K \sigma (x) + \G \mu(x) \right], 
	 	\end{equation}
		for all $x\in \Omega$ such that $u(x)\le \G(u^q d \sigma)(x)+ \G \mu(x)$.  In particular, 
		\eqref{upper-est-qm-cor-K} holds $d \sigma$-a.e.
	\end{cor}

	The following lemma provides bilateral pointwise estimates of solutions 
	to \eqref{sublin-sigma-mu} for  quasi-metric kernels, together with the existence criteria. 
	
		\begin{lemma}\label{exist-lemma} Let $\mu, \sigma \in \mathcal{M}^{+}(\Omega)$ ($\sigma\not=0$) and $0<q<1$.  Suppose $G$ is a  quasi-metric kernel 
		on $\Omega$. Then a nontrivial solution $u$ to \eqref{sublin-sigma-mu}  exists if and only if 
		 $\G \sigma <+\infty$, $\K \sigma <+\infty$, and $\G \mu < +\infty$ 
		 $d\sigma$-a.e., and satisfies  the bilateral pointwise estimates  
				\begin{align}
				\label{sublin-low-e} 
	& c \,   [ (\G \sigma(x))^{\frac{1}{1-q}} + 	 \K\sigma (x)]+  \G \mu(x) \le 	u (x),
	 \\
		& u(x)  \le  C \, [ (\G \sigma(x))^{\frac{1}{1-q}} + 	 \K\sigma (x) +  \G \mu(x) ],   
					 \label{sublin-up-e}  
					 \end{align}	
		$d \sigma$-a.e. in $\Omega$, 
		where $c, C$ are positive constants which depend only on $q$ and the 
		quasi-metric constant $\kappa$ of the kernel $G$. 	
	\end{lemma}

		%%%%%%%%%%%%%% 
		%%%%%%%%%%%%%%%%
Existence criteria can be stated in several equivalent forms using the following lemma.

	\begin{lemma}\label{equiv-lemma} Let $\mu, \sigma \in \mathcal{M}^{+}(\Omega)$ ($\sigma\not=0$) and $0<q<1$.  Suppose $G$ is a  quasi-metric kernel on $\Omega$. Then the following conditions 
	are equivalent:
	
(i) 	 $\G \sigma <+\infty$, $\K \sigma <+\infty$, and $\G \mu < +\infty$ 
		 $d\sigma$-a.e. 
		 
 (ii) $\G \sigma\not\equiv +\infty$, $\K \sigma \not\equiv +\infty$, and $\G \mu \not\equiv +\infty$. 
		 
(iii) Conditions \eqref{cond-a}--\eqref{cond-c} hold for some (or, equivalently, all) $x_0\in \Omega$ and 
$a>0$. 
	\end{lemma} 
	
		 \begin{remark} An alternative criterion for the existence of (super) solutions in the case of  quasi-metric kernels $G$ is deduced in Corollary \ref{equiv-cor} below. 
		 \end{remark}

	%%%%%%%%%%%%%%%%%%%%%%%

	\subsection{Bilateral bounds: quasi-metrically modifiable kernels}\label{q-m-m kernels} In this section, we give both lower estimates of super-solutions and upper 
	estimates of sb-solutions 
	to sublinear integral equations 
	\eqref{sublin-sigma-mu} with quasi-metrically modifiable kernels $G$ and  modifiers $m>0$, as defined in the Introduction.

	Examples of  quasi-metrically modifiable kernels can be found in \cite{An}, \cite{FNV}, \cite{FV2}, \cite{H}. In particular, for bounded domains $\Omega\subset\R^n$ satisfying the boundary Harnack principle, such as bounded Lipschitz, NTA or uniform domains, Green's kernels $G$ for the Laplacian and fractional Laplacian (in the case $0<\alpha\le 2$, $\alpha<n$) are quasi-metrically modifiable.

	Let $0<q<1$ and $\mu, \sigma \in \M(\Omega)$. 
	We discuss relations between solutions (as well as sub- and super-solutions) to the equations
	\begin{align}
	\label{int-eq-u}
			u & = \G(u^q d \sigma) + \G\mu, \quad 0< u<\infty  \quad d \sigma\textrm{-a.e.} \, \, \text{in} \, \, \Omega, 
				\\
		\label{int-eq-v}
			v & = \widetilde{\G} (v^q d \tilde{\sigma}) + \widetilde{\G} \tilde{\mu}, \quad 0< v<\infty   \quad d \tilde{\sigma}\textrm{-a.e.} \, \, \text{in} \, \, \Omega,
	\end{align}
	where $\widetilde{G}$ is the modified kernel \eqref{mod-ker} with modifier $m$,
	  and 
	  \begin{align}\label{v-u}
	  v \defeq \frac{u}{m}, \quad d \tilde{\sigma} \defeq m^{1+q} \, d\sigma, \quad d \tilde{\mu}= m \, d \mu.
	  \end{align}
	  Clearly, equations \eqref{int-eq-u} and \eqref{int-eq-v} are equivalent.
	  	
Similarly, the following two weighted norm inequalities are equivalent, 
	\begin{align}
		\label{weighted-G}
		\Vert \G \nu \Vert_{L^q(m \, d\sigma)} & \le \tilde{\varkappa}(\Omega) \, \int_\Omega m \, d\nu,  \qquad \forall \,  \nu \in \M(\Omega), 
	\\
		\label{weighted-G-t}
		\Vert \widetilde{\G} \nu \, \Vert_{L^q(\Omega, \tilde{\sigma})} & \le \tilde{\varkappa}(\Omega) \, \Vert \nu \Vert,  \qquad \forall \,  \nu \in \M_b(\Omega).
	\end{align}
	
	In \eqref{weighted-G}, without loss of generality we may assume 
	$ \int_\Omega m d\nu<\infty$.

\begin{remark} The least constant 
	$\tilde{\varkappa}=\tilde{\varkappa}(\Omega)$ is the same in \eqref{weighted-G} and \eqref{weighted-G-t}, 
	since the latter is 
	 an equivalent restatement of the former,  in terms of  $\widetilde{\G}$, $\tilde{\sigma}$ 
	in place of $\G$, $\sigma$.
\end{remark}

	Let  $\widetilde{B}=\widetilde{B}(x, r)$  be a quasi-metric ball in 
	$\Omega$ associated with the quasi-metric 
	$\tilde{d}=1/\widetilde{G}$, i.e., 		
	\begin{equation}\label{def-tilde-B} 
	\widetilde{B}(x, r) \defeq \left\{ y \in \Omega\colon \, \, \widetilde{G}(x, y)>1/r \right\}, \qquad  
	x \in \Omega, \, \, r>0.
		\end{equation}
		
		We denote by $\tilde{\varkappa}(\widetilde{B})$ the least constant in the localized  versions 
	of  inequalities  \eqref{weighted-G}, \eqref{weighted-G-t} with 
	$\sigma_{\widetilde{B}}$ in place of $\sigma$, and 
	$\tilde{\sigma}_{\widetilde{B}}$  in place of $\tilde{\sigma}$, respectively.

		Then the modified intrinsic potential $\widetilde{\K} \sigma$ is defined by 
		\begin{equation}\label{def-tilde-K} 
	\widetilde{\K} \sigma  (x) \defeq \int_0^\infty \frac{ [\tilde{\varkappa}(\widetilde{B}(x, r))]^{\frac{q}{1-q}} }{r^2} \, dr, \qquad x\in \Omega.
		\end{equation}
		
		The following lemma (\cite{Ver2}*{Lemma 5.1})   contains a lower bound for super-solutions to sublinear integral equations for quasi-metrically modifiable kernels  	$G$.

	\begin{lemma}
		\label{qm_lemma-lower}
		Let $G$ be a {\rm (QS)}   kernel on $\Omega$ with quasi-symmetry constant 
		$\mathfrak{a}$, 
		such that the modified kernel $\widetilde{G}$, defined by \eqref{mod-ker} with  modifier $m$,  
		satisfies the {\rm (WMP)} 
		with constant $\mathfrak{b}$.  Then any nontrivial super-solution $u$ to \eqref{int-eq-u}  satisfies 
		the estimate 
		\begin{equation}\label{lower-est-qm} 
	 		u\ge c \, \left(m \, \left[\frac{\G( m^q d \sigma)}{m}\right]^{\frac{1}{1-q} }
			+ \widetilde{\K} \sigma\right)     +  \, \G\mu\qquad d \sigma\textrm{-a.e.},
	 	\end{equation}
		where $c=c(q, \mathfrak{a}, \mathfrak{b})$ is a positive constant.
	\end{lemma}

	The following lemma (\cite{Ver2}*{Lemma 5.2})   provides a matching  upper  bound for sub-solutions to sublinear integral equations, for quasi-metrically modifiable kernels  	$G$.

	\begin{lemma}
		\label{qm_lemma-upper-mod}
		Let $G$ be a quasi-metrically modifiable  kernel on $\Omega$ with   modifier $m$.
			Then any sub-solution $u$ to \eqref{int-eq-u}   satisfies the estimate
		\begin{equation}\label{upper-est-qmm} 
	 		u\le C \, m \, \left( \left [ \frac{\G( m^q d \sigma)}{m}\right]^{\frac{1}{1-q} } 
			+ \widetilde{\K}  \sigma\right) + C \, \G\mu \quad d \sigma\textrm{-a.e.},
	 	\end{equation}
		where $C=C(q, \tilde{\kappa})$, and $\tilde{\kappa}$ is the quasi-metric constant   for the modified kernel 
		$ \widetilde{G}$.
	\end{lemma}

	We next consider the modifiers $m=g$ given by 
	\begin{equation}\label{g-mod-def}
	g(x) =g^{x_0} (x) \defeq \min\{1, G(x, x_0)\}, \quad x \in \Omega,
		\end{equation}
	where $x_0\in \Omega$ is a fixed pole.

	Let $G$ be a quasi-metric kernel  on $\Omega$, so that $d\defeq  1/G$ obeys 
	the quasi-triangle inequality \eqref{q-m-tr} 
	with quasi-metric with constant $\kappa$. 
	The proof of the following lemma (\cite{Ver2}*{Lemma 5.3})   is based on the so-called 
	 Ptolemy's inequality 
	  (see \cite{FNV}*{Sec. 3}):  for all $x_0, x, y, z \in \Omega$, 
		\begin{equation}\label{ptolemy} 
	 		d(x,y) d(x_0, z) \le 4 \kappa^2 \Big[d(x, z) d(y, x_0) + d(x_0, x) d(z, y)\Big]. 	\end{equation}

	\begin{lemma}\label{mod-lemma} Let  $G$ be a quasi-metric kernel on $\Omega$ with quasi-metric constant $\kappa$.  Let $x_0\in \Omega$, and let $g(x) = \min\{ 1, G(x, x_0)\}$. 
		Then 
		\begin{equation}\label{ker-mod-G}
		\widetilde{G}(x,y) = \frac{G(x,y)}{g(x) g(y)}
			\end{equation}  
		is a quasi-metric kernel on $\Omega$ 
	 with quasi-metric constant 	$4 \kappa^2$. In particular,  $\widetilde{G}$ 
satisfies the {\rm (WMP)} with constant $\mathfrak{b}=8 \kappa^2$. 
	\end{lemma}

	 In the next lemma, we give a criterion for the existence of (super) solutions 
	in the case of quasi-metrically modified kernels $G$.  
	
	\begin{lemma}\label{equiv-lemma-mod} Let $\mu, \sigma \in \mathcal{M}^{+}(\Omega)$  and $0<q<1$.  Suppose $G$ is a  quasi-metrically modifiable kernel on $\Omega$ with modifier  $m=g \in C(\Omega)$ defined by \eqref{g-mod-def}. Then there exists a nontrivial (super) solution to equation \eqref{int-eq-u} 
if and only if conditions \eqref{tilde-kappa} hold,  i.e., 
\begin{equation}\label{tilde-kappa-5}
\int_\Omega g \, d \mu<\infty \quad {\rm and} \quad  \tilde{\varkappa}(\Omega)<\infty,  
\end{equation}
where  $\tilde{\varkappa}(\Omega)$ is the least constant in the weighted norm inequality \eqref{weighted-G} with $m=g$. 
	\end{lemma} 
	
	We now go back to  quasi-metric kernels $G$. 
	The following corollary  is a direct consequence of 
	 Lemma \ref{equiv-lemma-mod}, since in this case 
	   the modified kernels $\widetilde{G}$ defined by 
	 \eqref{ker-mod-G} are also quasi-metric by Lemma \ref{mod-lemma}.

	\begin{cor}\label{equiv-cor} Let $\mu, \sigma \in \mathcal{M}^{+}(\Omega)$  and $0<q<1$.  Suppose $G$ is a  quasi-metric kernel on $\Omega$ such that $g \in C(\Omega)$, where  $g$ is defined by \eqref{g-mod-def}. Then there exists a nontrivial (super) solution to equation \eqref{int-eq-u} 
if and only if  \eqref{tilde-kappa-5} holds.
	\end{cor}

	%%%%%%%%%%%%%%%%%%%%
	%%%%%%%%%%%%%%%%%%%%

\section{Quasilinear equations with sub-natural growth terms}\label{Sec-6}

	\subsection{Nonlinear potential estimates} 	
	In this section,  we present  bilateral pointwise estimates of solutions 
 to   quasilinear elliptic equations of the type 
 \begin{equation}\label{inhom}
\begin{cases}
-\Delta_{p} u  = \sigma u^{q} + \mu, \quad 0<u<\infty \quad  d\sigma\textrm{-a.e.}  \quad \text{in} \;\; \R^n, \\
\liminf\limits_{x\rightarrow \infty}u(x) = 0,  
\end{cases}
\end{equation}
where $\mu$, $\sigma\in \M(\R^n)$, in the 
\textit{sub-natural growth} case $0<q<p-1$. 

Here all solutions $u$  are understood to be  $p$-superharmonic. This is a natural 
class of solutions to  \eqref{inhom}, since   $\nu\defeq \sigma u^{q} + \mu\ge 0$ in the sense of measures. 
We may  assume here, without loss of generality, that  
 $u\in L^{q}_{{\rm loc}}(\R^n, d\sigma)$, so that  
 $\nu \in \M(\R^n)$ (see \cite{Ver1}).

The notion of a $p$-superharmonic function is discussed below in a more general 
setting of $\mathcal{A}$-superharmonic functions associated with quasilinear  
equations involving the $\mathcal{A}$-Laplace operator, namely,  
 \begin{equation}\label{inhom-omega}
-{\rm div}(\mathcal{A}(x, {\rm D} u))=\nu \qquad {\rm in} \quad \Omega,
 \end{equation}
 where $\Omega \subseteq \R^n$,  $\nu \in \M(\Omega)$, 
 ${\rm D}$ is the generalized gradient defined below,  and $\mathcal{A}$ 
 obeys certain monotonicity and growth assumptions discussed below 
 (see  \cite{HKM}). 

\begin{remark} For equations \eqref{inhom-omega} with data $\omega \in \M(\Omega)$, 
the class of  $\mathcal{A}$-superharmonic solutions coincides in a sense with the class of \textit{local renormalized} solutions. We refer to \cite{KKT} for the proof of this important fact, and  the discussion of the literature on renormalized solutions. 
\end{remark}

 We will present matching upper and lower estimates of solutions to \eqref{inhom} 
 in terms of certain nonlinear potentials 
defined below. Our estimates hold for all $p$-superharmonic solutions $u$. In particular, they yield an existence criterion for solutions to  \eqref{inhom}. 

These results, obtained recently in \cite{Ver1}, are new even in the special case $\mu=0$, i.e., for the equation 
\begin{equation}\label{q-hom}
\begin{cases}
-\Delta_{p} u  = \sigma u^{q}, \quad u\ge 0  \quad \text{in} \;\; \R^n, \\
\liminf\limits_{x\rightarrow \infty}u(x) = 0.  
\end{cases}
\end{equation}

Equation \eqref{q-hom} was considered earlier in \cite{CV1}, but the upper pointwise estimate was obtained only for the 
\textit{minimal} solution $u$. Moreover, we will discuss   uniqueness results for \textit{reachable} solutions $u$ obtained very recently in \cite{PV2}.  

We will use the notion of the $p$-capacity for compact sets $K\subset \R^n$.

\begin{defn} Let $1<p<\infty$ and $K\subset \R^n$ be a compact set. 
The $p$-capacity of  $K$ is defined by 
\begin{equation}\label{p-capacity}
{\rm cap}_p (K)=\inf \left\{\int_{\R^n} |\nabla u|^p dx: \, \,  u \ge 1 \, \, \text{on} \, \, K, \quad u \in C^\infty_0(\R^n)\right \}.
\end{equation}
\end{defn}

Notice that the $p$-capacity on $\R^n$ is nontrivial only if $1<p<n$. 

We recall the following definition. 

\begin{defn} A measure $\sigma\in \M(\R^n)$ is said to be absolutely continuous with respect to the $p$-capacity  if 
$\sigma(K)=0$ whenever ${\rm cap}_p (K)=0$, for any compact set $K\subset\R^n$. In this case, we write $\sigma << {\rm cap}_p$. 
\end{defn}

We observe that the existence of a (super) solution  to 
\eqref{inhom} yields   $\sigma<< {\rm cap}_p$.   More precisely, it follows from \cite{CV1}*{Lemma 3.6} that 
if $u$ is a nontrivial super-solution  to \eqref{q-hom} in the case $0<q\le p-1$, then  
\begin{equation}\label{abs-cap0}
\sigma(K)\le \text{cap}_p(K)^{\frac{q}{p-1}}\left(\int_K u^q d\sigma\right)^{\frac{p-1-q}{p-1}},
\end{equation}
for all compact sets $K\subset \R^n$.

Among our main tools are certain \textit{nonlinear potentials} associated with \eqref{q-hom}. We recall that the Havin--Maz'ya--Wolff potential $\W_{\alpha, p}$  is defined, for  $1<p<\infty$ and 
$0<\alpha<\frac{n}{p}$, by \eqref{nonlin-MH}.
 
  In the special case
$\alpha=1$ ($1<p<n$) used in Sec. \ref{Sec-6}, this  nonlinear potential  
will be denoted by $\W_p$,  i.e., for $\mu \in \M(\R^n)$, we set 
\begin{equation}
	 \label{nonlin-MH-1}
	 \W_p\mu (x) \defeq \int_0^\infty \left[ \frac{\mu(B(x, \rho))}{\rho^{n-p}} \right]^{\frac{1}{p-1}}  \, \frac{d \rho}{\rho}, \qquad x \in \R^n, 
	 \end{equation}
	 where $B=B(x, \rho)$ 
	 is a Euclidean ball in $\R^n$ of radius $\rho$ centered at $x$.

	 For $\mu\in \M(\R^n)$, we consider the equation 
	 \begin{equation}\label{p-Lapl} 
\left\{ \begin{array}{ll}
- \Delta_p u = \mu, \quad u\geq 0  \quad \text{in } \R^n, \\
\displaystyle{\liminf_{x\rightarrow \infty}}\,  u = 0, 
\end{array}
\right.
\end{equation}

The following important global estimate, together with its local counterpart,  is due to T.~Kilpel\"{a}inen and J.~Mal\'y 
\cite{KM}. Suppose 
$u\ge 0$ is a $p$-superharmonic solution to \eqref{p-Lapl}. Then
\begin{equation}\label{k-m}
K^{-1}\W_p \mu (x) \leq u(x) \leq K \W_p \mu (x),
\end{equation} 
where $K=K(p, n)$ is a positive constant.

Moreover, it is known (see \cite{PV1}) that  a nontrivial solution $u$  to \eqref{p-Lapl} 
exists if and only if 
\begin{equation}\label{finiteness}
\int_{1}^{\infty}  \left[ \frac{\mu(B(0, \rho))}{\rho^{n-p}} \right]^{\frac{1}{p-1}} 
\frac{d \rho}{\rho}<\infty.
\end{equation} 
This is equivalent to $\W_p \mu (x)<\infty$ for some  $x \in \R^n$, or 
equivalently quasi-everywhere (q.e.) on $\R^n$  with respect to the $p$-capacity. 
In particular, \eqref{finiteness} may hold only in the case $1<p<n$, unless $\mu = 0$.

We next recall the definition the so-called \textit{intrinsic} nonlinear potential 
$\K_{p, q}$
associated with \eqref{q-hom}, which was introduced in \cite{CV1}.  

To define $\K_{p, q} \sigma$ for $\sigma \in \M(\R^n)$, 
we first consider the weighted norm inequality 
\begin{equation} \label{weight-lap-glob}
\left(\int_{\R^n} |\varphi|^q \, d \sigma\right)^{\frac 1 q} \le C \,  \Vert \Delta_p \varphi\Vert^{\frac{1}{p-1}},  
\end{equation} 
for all test functions $\varphi$  which are $p$-superharmonic in $\R^n$,  
and  such that $\displaystyle{\liminf_{x \rightarrow  \infty}} \, \varphi(x)=0$. Here 
$-\Delta_p \varphi =\nu$ is the Riesz  measure of $\varphi$, and 
without loss of generality we may assume  $\Vert \Delta_p \varphi\Vert 
=\nu(\R^n)<\infty$. 

By $\varkappa(\R^n)$ we denote the least constant in \eqref{weight-lap-glob}.

The  nonlinear potential $\K_{p, q} \sigma$ is defined in terms of the  localized version of \eqref{weight-lap-glob},    namely, 
\begin{equation} \label{weight-lap}
\left(\int_{B} |\varphi|^q \, d \sigma\right)^{\frac 1 q} \le C \,  \Vert \Delta_p \varphi\Vert^{\frac{1}{p-1}},  
\end{equation} 
where $B$ is a Euclidean ball in $\R^n$, for the same class of $p$-superharmonic  test functions $\varphi$ in $\R^n$  
as in \eqref{weight-lap-glob}.

By $\varkappa(B)$ we denote the least constant in \eqref{weight-lap} 
with $\sigma_B$ in place of $\sigma$, where $\sigma_B=\sigma\vert_B$ is restricted to a ball 
$B=B(x, \rho)$.  
Then the \textit{intrinsic} nonlinear potential $\K_{p, q}\sigma$ is defined by 
\begin{equation} \label{potentialK}
\K_{p, q}  \sigma (x)  \defeq  \int_0^{\infty} \left[\frac{ \varkappa(B(x, \rho))^{\frac{q(p-1)}{p-1-q}}}{\rho^{n- p}}\right]^{\frac{1}{p-1}}\frac{d\rho}{\rho}, \quad x \in \R^n.
\end{equation} 

As was noticed in \cite{CV1},  $\K_{p, q}  \sigma \not\equiv + \infty$ if and only if 
\begin{equation}\label{suffcond1}
\int_1^{\infty} \left[\frac{ \varkappa(B(0, \rho))^{\frac{q(p-1)}{p-1-q}}}{\rho^{n- p}}\right]^{\frac{1}{p-1}}\frac{d\rho}{\rho} < \infty.  
\end{equation}
If \eqref{suffcond1} holds, then actually $\K_{p, q}  \sigma < + \infty$ $d \sigma$-a.e., and q.e. with respect to the $p$-capacity. 

\begin{remark} In the case 
$p=2$ and $\Omega = \R^n$, the potential $\K_{2, q} $ is closely related to the 
 nonlinear potential $\K$,   defined   by \eqref{nonlin_pot}.  
Notice that in \eqref{nonlin_pot} in this special case, $B(x, r)$ stands 
for a quasi-metric ball with respect to the quasi-metric $d(x, y)=|x-y|^{n-2}$, 
 $n\ge 3$, whereas 
in \eqref{potentialK}, $B(x, \rho)$ is a Euclidean ball. Hence, using the substitution 
$\rho=r^{\frac{1}{n-2}}$, we see that  $\K=(n-2) \K_{2, q}$. 
\end{remark}

We now state the main result of \cite{Ver1}, which establishes global bilateral estimates and existence criteria for all solutions to \eqref{inhom}. Besides 
 \eqref{finiteness}, \eqref{suffcond1}, we will use additionally  the condition  $\W_p \sigma \not\equiv + \infty$, i.e.,
\begin{equation}\label{suffcond}
\int_{1}^{\infty}  \left[\frac{\sigma(B(0, \rho))}{\rho^{n-p}}\right]^{\frac{1}{p-1}}
 \frac{d \rho}{\rho} <\infty.  
\end{equation}

\begin{theorem}\label{thm:main1}
Let $1<p<n$, $0<q<p-1$, and $\mu, \sigma \in M^{+}(\R^n)$. There exists a 
nontrivial $p$-superharmonic  solution $u$ to \eqref{inhom} if and only if conditions \eqref{finiteness}, 
\eqref{suffcond1}, and \eqref{suffcond} hold. Moreover, any such a solution $u$  satisfies the estimates 
\begin{align}\label{main-b1} 
 C_1 \, & \left[  (\W_p \sigma(x))^{\frac{p-1}{p-1-q}} + \K_{p, q} \sigma(x)  + \W_p \mu (x)\right] \le u(x) \\ 
 u(x) &  \le C_2 \, \left[  (\W_p \sigma(x))^{\frac{p-1}{p-1-q}} + \K_{p, q} \sigma(x) + \W_p \mu (x)\right],  \label{main-b2} 
\end{align}
$d \sigma$-a.e. on $\R^n$, with positive constants $C_i=C_i(p, q, n)$   ($i=1, 2$). 

If $n \leq p < \infty$, there are no nontrivial $p$-superharmonic solutions.  
\end{theorem}

\begin{remark}\label{p-sup-sub} As in the case of sublinear problems discussed above,  \eqref{main-b1}  holds for all nontrivial 
$p$-superharmonic super-solutions, whereas  \eqref{main-b2} holds for all $p$-superharmonic sub-solutions.
\end{remark}

\subsection{Equations involving $\mathcal{A}$-Laplace operators}

Our next goal is to introduce  the notion of a \textit{reachable} solution to equation \eqref{p-Lapl}, 
and discuss criteria of existence and uniqueness for reachable solutions to 
equation \eqref{inhom} in the  case $0<q<p-1$. 

We actually consider more general quasilinear  $\mathcal{A}$-Laplace operators 
in place of $\Delta_p$.

Let   $\mathcal{A}\colon  \R^n \times \R^n\rightarrow\R^n$ be a Carath\'eodory function such that 
the map $x \rightarrow \mathcal{A}(x,\xi)$  is  measurable~ for~ all~ $\xi\in\R^n,$ and 
the~ map ~ $\xi\rightarrow \mathcal{A}(x,\xi)$  ~is~ continuous~ for~ a.e. $x\in\R^n$.
We also assume that   there are constants $0<\alpha\leq\beta<\infty$ and $1<p<n$ such that for a.e. $x$ in $\R^n$,
\begin{equation}\label{structure}
\begin{aligned}
& \mathcal{A}(x,\xi)\cdot\xi\geq \alpha |\xi|^p,\quad |\mathcal{A}(x,\xi)|\leq \beta| \xi|^{p-1},  \quad \forall\, \xi 
\in \R^n, \\
& [\mathcal{A}(x,\xi_{1})-\mathcal{A}(x,\xi_{2})]\cdot(\xi_{1}-\xi_{2})>0, \quad \forall\, \xi_{1}, \xi_2 \in \R^n, \, \, 
\xi_{1}\not = \xi_{2}.
\end{aligned}
\end{equation}

For  the uniqueness results, we will assume additionally   
\begin{equation}\label{homogeneity}
\mathcal{A}(x,\lambda \xi)=\lambda^{p-1} \mathcal{A}(x, \xi), \qquad \forall \, \xi 
\in \R^n,  \,
\lambda>0.
\end{equation}
Condition \eqref{homogeneity} is often used in the literature ( \cite{HKM}, \cite{KM}). 

The special case $\mathcal{A}(x,\xi)=|\xi|^{p-2}\xi$ gives  the $p$-Laplacian 
$\Delta_p$. 

 For an open set $\Omega\subset\R^n$, it is well known that every weak solution $u\in W^{1,\,p}_{{\rm loc}}(\Omega)$  to the
equation
\begin{eqnarray}
\label{homo}
-\text{div}\mathcal{A}(x,\nabla u)=0 \qquad \textrm{in} \, \, \Om
\end{eqnarray}
has a continuous representative. Such continuous solutions are said to be
$\mathcal{A}$-$harmonic$ in $\Omega$. If $u\in W_{{\rm loc}}^{1,\,p}(\Omega)$ and
\begin{eqnarray*}
	\int_{\Om}\mathcal{A}(x,\nabla u)\cdot\nabla\varphi \, dx\geq 0,
\end{eqnarray*}
for all nonnegative $\varphi\in C^{\infty}_{0}(\Om)$, i.e., $-{\rm div}\mathcal{A}(x,\nabla u)
\geq 0$ in the distributional sense, then $u$ is called a {\it super-solution}
to (\ref{homo}) in $\Om$.\\
\indent A  function $u\colon\Om\rightarrow (-\infty, \infty]$ is called
$\mathcal{A}$-$superharmonic$ if $u$ is not identically infinite in each connected component
of $\Om$, $u$ is lower semicontinuous, and for all open sets  $D$ such that
${\overline D}\subset\Om$, and all functions $h\in C(\overline{D})$, $\mathcal{A}$-harmonic in $D$, it follows that
$h\leq u$ on $\partial D$ implies $h\leq u$ in $D$.

It is well known  that if  $u$ is an $\mathcal{A}$-superharmonic function, then 
for any
  $k>0$,   
  its truncation $u_k=\min\{u,k\}$ is $\mathcal{A}$-superharmonic as well. Moreover,   $u_k\in W^{1,\,p}_{{\rm loc}}(\Om)$ (see \cite{HKM}). 
We will need the notion of the  weak (generalized) gradient of $u$  defined by
\begin{eqnarray*}
	{\rm D} u \defeq \lim_{k\rightarrow\infty} \, \nabla  [ \, \min\{u,k\}] \qquad {\rm a.e.}  \,\, \textrm{in} \, \, \Om.
\end{eqnarray*}

We observe that ${\rm D} u$ gives the usual distributional gradient $\nabla u$
if either $u\in L_{{\rm loc}}^{\infty}(\Om)$ or $u\in W^{1,\,1}_{{\rm loc}}(\Om)$. Moreover, there exists a unique measure $\mu=\mu[u]\in \M(\Om)$ called the \textit{Riesz measure} of $u$ such that
\begin{equation}\label{eq-mu}
	-{\rm div}\mathcal{A}(x, {\rm D} u)=\mu \qquad \textrm{in} \, \, \Om.
\end{equation}

Let $\mu \in \M(\R^n)$. 
We first treat the problems of  existence and uniqueness  of $\mathcal{A}$-superharmonic solutions  to the equation
\begin{equation}\label{Basic-PDE}
\left\{ \begin{array}{ll}
-{\rm div}\, \mathcal{A}(x, {\rm D} u) = \mu, \quad u\geq 0  \quad \text{in } \R^n, \\
\displaystyle{\liminf_{x\rightarrow \infty}}\,  u = 0, 
\end{array}
\right.
\end{equation}
where $\mu \in \M(\R^n)$.

We notice that, under conditions \eqref{structure}, the  Kilpel\"{a}inen--Mal\'y  estimates \eqref{k-m}  hold 
for  all $\mathcal{A}$-superharmonic solutions to \eqref{Basic-PDE}, i.e., 
\begin{equation}\label{k-m-a}
K^{-1}\W_p \mu (x) \leq u(x) \leq K \W_p \mu (x),
\end{equation} 
where   $K=K(p, n, \alpha, \beta)$ is a positive constant (see \cite{KM}, \cite{KuMi}). 

Moreover, it is known \cite{PV1} that  a necessary and sufficient condition for  
the existence of an $\mathcal{A}$-superharmonic  solution to 
\eqref{Basic-PDE}  is given by \eqref{finiteness}, as in the case of equation 
\eqref{p-Lapl} 
for the $p$-Laplacian. 

We next observe 
that Theorem \ref{thm:main1} has a complete analogue (see \cite{Ver1}*{Remark 4.3(2)} 
for the equation 
\begin{equation}\label{sub-quasi-pde-a}
\begin{cases}
-{\rm div}  \mathcal{A} (x, {\rm D} u)  = \sigma u^{q} + \mu, \quad 0<u<\infty \, \, d\sigma\textrm{-a.e.}  \quad \text{in} \;\; {\mathbb R}^n, \\
\liminf\limits_{x\rightarrow \infty}u(x) = 0,  
\end{cases}
\end{equation}
where $\mu$, $\sigma\in \mathcal{M}^{+}({\mathbb R}^n)$.

	A natural existence criterion for equation \eqref{sub-quasi-pde-a}, together 
 with bilateral pointwise estimates, 
	is contained in the following theorem.

	\begin{remark} Some corrections in the proof of existence, provided in \cite{PV2}*{Lemma 4.1}, are 
	needed in the comparison principle  for 
	$\mathcal{A}$-superharmonic functions in 
	bounded domains $\Omega$  (\cite{CV1}*{Lemma 5.2}) used in 
	the constructions of solutions in \cite{CV1} for $\mu=0$, and \cite{Ver1} 
	for $\mu\not=0$.
	\end{remark}

\begin{theorem}\label{thm:main1a}
Let $1<p<n$, $0<q<p-1$, and $\mu, \sigma \in M^{+}({\mathbb R}^n)$. Suppose that 
 $\mathcal{A}$ satisfies conditions \eqref{structure}. Then there exists a 
nontrivial  $\mathcal{A}$-superharmonic solution $u$ to \eqref{sub-quasi-pde-a} if and only if conditions \eqref{finiteness}, 
\eqref{suffcond1}, and \eqref{suffcond} hold. 

Moreover, any nontrivial solution $u$  satisfies estimates \eqref{main-b1} 
with positive constants $C_i=C_i(p, q, n, \alpha, \beta)$ $(i=1, 2)$. 

If $n \leq p < \infty$, then there exist no nontrivial $\mathcal{A}$-superharmonic  solutions to \eqref{sub-quasi-pde-a}.  
\end{theorem}

An analogue of Remark \ref{p-sup-sub} remains true for $\mathcal{A}$-superharmonic  sub- and super-solutions 
of \eqref{sub-quasi-pde-a}.

\subsection{Reachable solutions to basic quasilinear equations}

In this section, we will define \textit{reachable} $\mathcal{A}$-superharmonic solutions to \eqref{Basic-PDE} for which existence is obtained under the sole condition 
\eqref{finiteness}, and uniqueness is ensured if additionally 
$\mu<<{\rm cap}_p$.    We recall that \eqref{finiteness} is necessary for 
the existence of an $\mathcal{A}$-superharmonic solution.

We first notice the existence of a minimal solution to  \eqref{Basic-PDE} if 
\eqref{finiteness} holds and  $\mu<<{\rm cap}_p$.
\begin{theorem}\label{miniexist}
	Let  $\mu \in \M(\R^n)$, where $\mu<< {\rm cap}_p$. Suppose that 
	\eqref{finiteness} holds.
	Then there exists a  minimal $\mathcal{A}$-superharmonic solution to equation \eqref{Basic-PDE}.
\end{theorem}

\begin{remark} 1. It is not known if condition \eqref{finiteness} alone is enough for the existence of the minimal solution in Theorem \ref{miniexist}. 

2. It is also not known whether, under  condition \eqref{finiteness} combined with $\mu << {\rm cap}_p$,  an $\mathcal{A}$-superharmonic solution to \eqref{Basic-PDE} is unique, and hence coincides with the minimal solution. Some partial results in this direction will be discussed below.
 \end{remark}

We now consider the following notion of  a reachable solution (see \cite{DMM}*{Definition 2.3} in the case of bounded domains) suitable for our purposes.

\begin{defn}\label{reach}
	Let  $\mu \in \M(\R^n)$.
	A function $u\colon \R^n \rightarrow [0, +\infty]$  is said to be an  $\mathcal{A}$-superharmonic reachable solution to equation \eqref{Basic-PDE} 
	if $u$ is an	$\mathcal{A}$-superharmonic solution of \eqref{Basic-PDE}, and there exist two sequences $\{u_i\}$ and $\{\mu_i\}$, $i=1,2, \dots$,  such that 
	
	\noindent ${\rm (i)}$ Each $\mu_i\in \M(\R^n)$ is  compactly  supported  in $\R^n$, and  $\mu_i \leq  \sigma $; 
	%and $\sigma_i \rightarrow \sigma$ weakly as measures in $\R^n$.  
	
	\noindent ${\rm (ii)}$ Each $u_i$ is 	an $\mathcal{A}$-superharmonic solution of \eqref{Basic-PDE} with  $\mu_i$ in place of $\mu$;

	\noindent ${\rm (iii)}$ $u_i\rightarrow u$ a.e. in $\R^n$.
\end{defn}

\begin{remark}  The requirement that $\mu_i \leq \mu$ in Definition \ref{reach} is important in the proof of uniqueness  in the case $\mu<<{\rm cap}_p$. 
 \end{remark}
 
 The next theorem is the main result of \cite{PV2} on reachable solutions to equation \eqref{Basic-PDE}.

\begin{theorem}\label{EandU}
	Suppose $\mu\in \M(\R^n)$, and 
	\eqref{finiteness} holds. Then there exists an  $\mathcal{A}$-superharmonic reachable solution to \eqref{Basic-PDE}. Moreover, if  additionally $\mu<< {\rm cap}_p$, then any 
	$\mathcal{A}$-superharmonic reachable solution is unique and coincides with the minimal solution.  
\end{theorem}

\begin{remark}\label{rem 5-10} 
For $\mu\in \M_b(\R^n)$, it is known \cite{PV2}*{Theorem 3.12} that  any  $\mathcal{A}$-superharmonic solution $u$ to  \eqref{Basic-PDE} is unique, and 
coincides with the minimal  $\mathcal{A}$-superharmonic solution.
\end{remark}

The next  theorem proved in \cite{PV2}   
shows that all $\mathcal{A}$-superharmonic solutions to  \eqref{Basic-PDE} are  reachable, provided the condition  $\displaystyle{\liminf_{x\rightarrow \infty}}\,  u = 0$  is replaced with 
$\displaystyle{\lim_{x\rightarrow \infty}}\,  u = 0$.

\begin{theorem}\label{thm-lim} Let 
 $\mu \in \M(\R^n)$, and 
$\mu<<{\rm cap}_p$.  	Suppose that  $u$ is an $\mathcal{A}$-superharmonic solution  to the equation
	\begin{equation}\label{Basic-PDE2}
	\left\{ \begin{array}{ll}
	-{\rm div}\, \mathcal{A}(x, {\rm D} u) = \mu, \quad u\geq 0  \quad \text{in } \R^n, \\
	\displaystyle{\lim_{x\rightarrow \infty}}\,  u = 0, 
	\end{array}
	\right.
	\end{equation} 
	Then $u$ is the unique $\mathcal{A}$-superharmonic solution of \eqref{Basic-PDE2}, which coincides with the minimal 
$\mathcal{A}$-superharmonic reachable solution of \eqref{Basic-PDE}.
\end{theorem}

\subsection{Existence of reachable solutions: 
 sub-natural growth} We now discuss existence of nontrivial \textit{reachable}  $\mathcal{A}$-superharmonic solutions $u$ to equation \eqref{sub-quasi-pde-a}.  
	As above, we assume without loss 
of generality that $u \in L^q_{{\rm loc}} ({\mathbb R}^n, \sigma)$, so that 
	 $\sigma u^q + \mu \in \mathcal{M}^{+}({\mathbb R}^n)$, and $\sigma << {\rm cap}_p$.

We first consider homogeneous equations \eqref{sub-quasi-pde-a} with $\mu=0$. 
By  Theorem \ref{thm:main1a}, there exists a nontrivial $\mathcal{A}$-superharmonic solution $u$ if and only if $\mathbf{W}_{p}  \sigma \not\equiv + \infty$, $\mathbf{K}_{p, q}  \sigma \not\equiv + \infty$, i.e., 
conditions \eqref{finiteness}, \eqref{suffcond1}   hold. This theorem is complemented  by the following statement proved in \cite{PV2}*{Theorem 4.2}. 

\begin{theorem}\label{exist} Let $0<q<p-1$, and let $\sigma\in \mathcal{M}^{+}({\mathbb R}^n)$.  
	Then the nontrivial minimal $\mathcal{A}$-superharmonic solution $u$ of 
	\begin{equation}\label{sub-quasi-pde}
	\left\{ \begin{array}{ll}
	-{\rm div}\, \mathcal{A}(x, {\rm D} u) = \sigma u^q, \quad u\ge 0   \quad \text{in } {\mathbb R}^n, \\
	\displaystyle{\liminf_{x\rightarrow \infty}}\,  u = 0, 
	\end{array}
	\right.
	\end{equation}
 constructed in the proof of \cite{CV1}*{Theorem 1.1} under  conditions  \eqref{finiteness}, \eqref{suffcond1}, is in fact an $\mathcal{A}$-superharmonic reachable solution. 
\end{theorem}

The following statement (\cite{PV2}*{Theorem 4.3}) provides an analogue 
of Theorem \ref{exist} for reachable solutions to inhomogeneous equations 
\eqref{sub-quasi-pde-a}. 
In contrast  to the construction of an $\mathcal{A}$-superharmonic  solution 
(not necessarily reachable) 
in \cite{Ver1}*{Theorem 1.1}, the proof  for $\mu\not=0$ is different, and  relies on the extra 
assumption $\mu<<{\rm cap}_p$.  

\begin{theorem} Let $0<q<p-1$, and let $\mu, \sigma\in \mathcal{M}^{+}({\mathbb R}^n)$, where 
$\mu<<{\rm cap}_p$.  
	Then, under conditions  \eqref{finiteness}, 
\eqref{suffcond1}, and \eqref{suffcond}, 
	there exists a  nontrivial minimal reachable $\mathcal{A}$-superharmonic solution of \eqref{sub-quasi-pde-a}.  
\end{theorem}

\subsection{Uniqueness of reachable solutions: 
 sub-natural growth} In conclusion, we discuss the uniqueness property for (reachable) solutions of \eqref{sub-quasi-pde-a}. The  following main theorem was obtained \cite{PV2}*{Theorem 4.4}. 

\begin{theorem}\label{main-thm-q-p} Let $0<q<p-1$, and let $\mu, \sigma\in \M(\R^n)$, where $\mu<<{\rm cap}_p$. Suppose $\mathcal{A}$ satisfies conditions \eqref{structure} and \eqref{homogeneity}. 
Then nontrivial $\mathcal{A}$-superharmonic reachable solutions 
of \eqref{sub-quasi-pde-a} are unique. 
\end{theorem}

It is known that in some cases listed below the restriction to reachable solutions in this uniqueness property can 
be dropped. 

\begin{remark} 1.  In the case $p=2$ of  Theorem \ref{main-thm-q-p},  for linear uniformly elliptic operators  $\mathcal{L}$ 
with bounded measurable coefficients given by 
\eqref{typ-L},   all $\mathcal{L}$-superharmonic solutions  of \eqref{sub-quasi-pde-a} are  unique,  without the extra restriction $\mu<<{\rm cap}$.

2.   All nontrivial $\mathcal{A}$-superharmonic solutions in Theorem \ref{main-thm-q-p} 
 are unique if any one of the following conditions hold (\cite{PV2}*{Corollary 4.5}):

\noindent {\rm (i)} 	$\displaystyle{\lim_{x \rightarrow \infty}} u(x) =0$;

\noindent {\rm (ii)}   $u \in L^q (\R^n, d \sigma)$  and $ \mu \in \M_b(\R^n)$;
			
		\noindent {\rm (iii)} $|\nabla u| \in L^{p}(\R^n)$.
		\end{remark}

		%%%%%%%%%%%%%Bibliography%%%%%%%%%%%%%
%%%%%%%%%%%%%%%%%%%%%%%%%%%%%%%%%%%%%%%%%%%%%%%%


\begin{thebibliography}{99}

  
 \bibitem{AH}   {\sc D. R. Adams and L. I. Hedberg}, {\em Function Spaces and Potential Theory},
  Grundlehren der math. Wissenschaften {\bf 314}, Berlin--Heidelberg--New York, Springer, 1996.
    
\bibitem{Ai} {\sc H. Aikawa}, {\em Boundary Harnack principle and Martin boundary for a
uniform domain,} J. Math. Soc. Japan {\bf 53} (2001) 119--145.
  
\bibitem{An} {\sc A. Ancona}, {\em Some results and examples about
the behaviour of harmonic functions and Green's functions with
respect to second order elliptic operators,} Nagoya Math. J. {\bf 165} (2002) 123--158.

 \bibitem{BoOr} {\sc L. Boccardo and L. Orsina}, {\em Sublinear equations in {$L^s$}},  
Houston J. Math. \textbf{20} (1994), 
99--114. 

\bibitem{Brelot}{\sc M.Brelot}, {\em Lectures on Potential Theory},   Lectures on Math. {\bf 19}, Tata
  Institute, Bombay, 1960.

  
 \bibitem{BK}
 {\sc H.~Brezis and   S. Kamin}, {\em Sublinear elliptic equations on $\mathbb{R}^n$},  
 Manuscr.  Math. {\bf 74} (1992),  87--106. 
 
 \bibitem{BO}
{\sc H.~Brezis and L.~Oswald}, {\em Remarks on sublinear elliptic equations}, J. Nonlin. Anal.   
TMA {\bf 10} (1986), 55--64.
  
 \bibitem{CV1} {\sc Dat T. Cao and I. E. Verbitsky}, {\em Nonlinear elliptic equations and intrinsic potentials of Wolff type,} J. Funct. Anal.  {\bf 272} (2017), 112--165. 
 
 \bibitem{DMM} {\sc  G. Dal Maso  and A.  Malusa},    
{\em  Some properties of reachable solutions of nonlinear elliptic equations with measure data,}  Ann. Scuola Norm. Super. Pisa,   Ser. IV, {\bf 25} (1997),  375--396. 

\bibitem{FNV}{\sc M.~Frazier, F.~Nazarov, and I.~Verbitsky}, {\em Global estimates for
  kernels of Neumann series and Green's functions},  J. London Math. Soc. {\bf 90}   (2014),  903--918.

\bibitem{FV1}
{\sc M.~Frazier and I.~Verbitsky}, {\em Postive solutions to Schr\"odinger's equation and the exponential integrability of the balayage}, Ann. Inst. Fourier (Grenoble) {\bf 67} (2017), 1393--1425.

\bibitem{FV2}
{\sc M.~Frazier and I.~Verbitsky}, {\em Positive solutions and harmonic measure for Schr\"{o}dinger operators in uniform domains}, Pure Appl. Funct. Anal. (Yokohama Publ.), {\bf 7} (2022), 993--1023. 
 
\bibitem{Fug}{\sc B.~Fuglede}, {\em On the theory of potentials in locally compact spaces}, Acta Math. {\bf 103} (1960), 139--215.

\bibitem{G}{\sc E.~Gagliardo}, {\em On integral transformations with positive kernel},
  Proc. Amer. Math. Soc. {\bf 16} (1965), 429--434.
  
 \bibitem{GSV}  {\sc A.~Grigor'yan,  Y.~Sun, and I.~Verbitsky,} {\em Superlinear elliptic inequalities on manifolds,} J. Funct.  Anal. {\bf 278} (2020),  no. 9, Art. 108444.
    
\bibitem{GV} {\sc A.~Grigor'yan and I.~Verbitsky,} {\em Pointwise estimates of solutions to nonlinear equations for nonlocal operators,} Ann. Scuola Norm. Super. Pisa, Ser. V, {\bf 20} (2020),  721--750. 

\bibitem{H} {\sc W.~Hansen}, {\em Uniform boundary Harnack principle and generalized
  triangle property}, J. Funct. Anal. {\bf 226} (2005), 452--484.

    
 \bibitem{HW} {\sc L. I. Hedberg and T. H. Wolff}, {\em Thin sets in nonlinear potential theory},  
Ann. Inst. Fourier (Grenoble) {\bf 33} (1983), 161--187.

\bibitem{HKM} {\sc J. Heinonen, T. Kilpel\"ainen,  and   O. Martio},  
{\em  Nonlinear Potential Theory of Degenerate Elliptic Equations}, 
 Oxford Univ. Press, Oxford, 1993. 

  
\bibitem{KV} {\sc N.~J. Kalton and I.~E. Verbitsky}, {\em Nonlinear equations and weighted norm inequalities}, Trans. Amer. Math. Soc. {\bf 351} (1999), 3441--3497.

\bibitem{KKT} {\sc T. Kilpel\"ainen, T. Kuusi, and A. Tuhola-Kujanp\"a\"a},   
{\em  Superharmonic functions are locally renormalized solutions}, 
Ann. Inst. H. Poincar\'e Anal. Non Lin\'eaire {\bf 28} (2011), 775--795.


\bibitem{KM} {\sc T. Kilpel\"ainen and J. Mal\'y},   
{\em  The Wiener test and potential estimates for quasilinear elliptic equations},
Acta Math. {\bf 172} (1994), 137--161.	


\bibitem{Kr} {\sc  M.~A. Krasnoselskii}, {\em Positive Solutions of Operator Equations},  P. Noordhoff, Groningen,  1964. 

\bibitem{KuMi} {\sc T. Kuusi and G. Mingione}, {\em Guide to nonlinear potential estimates},  
Bull. Math. Sci. {\bf 4} (2014), 
1--82.


 
  \bibitem{MV} {\sc M. Marcus and L. V\'eron,} {\em Nonlinear Second Order Elliptic Equations Involving Measures}, Walter de Gruyter, Berlin--Boston, 2014. 

 \bibitem{Maz} {\sc V. Maz'ya}, {\em Sobolev Spaces, with Applications to Elliptic Partial Differential Equations}, 2nd Augm. Ed., Grundlehren der math. Wissenschaften \textbf{342}, Springer, Berlin, 2011.
  
 \bibitem{MH} {\sc V.~G. Maz'ya and V.~P. Havin}, {\em Non-linear potential theory},  Uspekhi Mat.
Nauk {\bf 27} (1972), 67--138; English transl.  Russ. Math. Surveys {\bf 27} (1972), 71--148.

 
 \bibitem{PV1}{\sc N.~C. Phuc and I.~E. Verbitsky}, 
  {\em Quasilinear and Hessian equations of Lane--Emden type}, 
  Ann. Math. {\bf 168} (2008), 859--914.
 
 
 \bibitem{PV2}{\sc N.~C. Phuc and I.~E. Verbitsky}, 
 {\em Uniqueness of entire solutions to quasilinear equations of $p$-Laplace type}, 
 arXiv:2208.13272.
 

\bibitem{QV2}{\sc S.~Quinn and I.~E. Verbitsky}, {\em A sublinear version of Schur's lemma and elliptic PDE}, Analysis \& PDE   {\bf 11}  (2018) 439--466.


\bibitem{SV}{\sc A.~Seesanea and I.~E. Verbitsky}, {\em Finite energy solutions to inhomogeneous nonlinear elliptic equations with sub-natural growth terms}, Adv. Calc. Var. {\bf 13} (2020)  53--74. 

\bibitem{V1}{\sc I.~E. Verbitsky}, {\em Sublinear equations and Schur's test 
for integral operators},  50 Years with Hardy Spaces, a Tribute to Victor Havin, eds. A. Baranov et al.,  Birkh\"{a}user, Operator Theory: Adv. Appl. {\bf 261}  (2018) 465--482.



\bibitem{Ver1}{\sc I.~E. Verbitsky},  {\em Bilateral estimates of solutions to  quasilinear elliptic equations with sub-natural growth terms}, Adv. Calc. Var.  (published online),  DOI: 10.1515/acv-2021-0004,  arXiv:2101.02368. 

 \bibitem{Ver2} {\sc I.~E. Verbitsky}, {\em 
Global pointwise estimates of positive solutions to sublinear equations},  
  St. Petersburg Math. J. (to appear), reprint of  Algebra i Analiz {\bf 34}, no. 3 (2022), 297--330,    
 arXiv:2203.0253.  

\end{thebibliography}
\end{document}